%
%
%
\input amstex.tex
\documentstyle{amsppt}
\magnification=1200
\baselineskip=13pt
\hsize=6.5truein
\vsize=8.9truein
\parindent=20pt
\countdef\sectionno=1
\countdef\eqnumber=10
\countdef\theoremno=11
\countdef\countrefno=12
\countdef\cntsubsecno=13
\sectionno=0

\def\newsection{\global\advance\sectionno by 1
                \global\eqnumber=1
                \global\theoremno=1
                \global\cntsubsecno=0
                \the\sectionno.\ }

\def\newsubsection{\global\advance\cntsubsecno by 1
                   \the\sectionno.\the\cntsubsecno.\ }

\def\theoremname#1{\the\sectionno.\the\theoremno
                   \xdef#1{{\the\sectionno.\the\theoremno}}
                   \global\advance\theoremno by 1}

\def\eqname#1{\the\sectionno.\the\eqnumber
              \xdef#1{{\the\sectionno.\the\eqnumber}}
              \global\advance\eqnumber by 1}

\global\countrefno=1

\def\refno#1{\xdef#1{{\the\countrefno}}
\global\advance\countrefno by 1}

\def\thmref#1{#1}

\def\R{{\Bbb R}}

\def\C{{\Bbb C}}
\def\Z{{\Bbb Z}}
\def\T{{\Bbb T}}
\def\Zp{{\Bbb Z}_+}

\def\Asu{A_q(SU(2))}
\def\U{U_q({\frak{sl}}(2,\C))}
\def\Usu{U_q({\frak{su}}(2))}
\def\Unc{U_q({\frak{su}}(1,1))}

\def\al{\alpha}
\def\be{\beta}
\def\ga{\gamma}
\def\de{\delta}
\def\De{\Delta}

\def\th{\theta}
\def\la{\lambda}

\def\vp{\varphi}
\def\Hi{\ell^2(\Zp)}

\refno{\AlSa}
\refno{\AlSaC}
\refno{\AskeW}
\refno{\Bere}
\refno{\BurbK}
\refno{\CharP}
\refno{\DunkPJM}
\refno{\DunkCJM}
\refno{\HTFtwee} 
\refno{\GaspR}
\refno{\GranZ}
\refno{\GrozK}
\refno{\KalnMM}
\refno{\KlimK}
\refno{\KoekS}
\refno{\KoelPAMS}
\refno{\KoelAAM}
\refno{\KoorJMP}
\refno{\KoorProb}
\refno{\KoorContM}
\refno{\KoorZSE}
\refno{\MassR}
\refno{\VdJeug}
\refno{\VileK}
\topmatter
\title Convolutions for orthogonal polynomials from
Lie and quantum algebra representations\endtitle
\author H.T. Koelink and J. Van der Jeugt\endauthor
\rightheadtext{Convolutions for orthogonal polynomials}
\affil Report 96-11, Universiteit van Amsterdam\endaffil
\address Vakgroep Wiskunde, Universiteit van Amsterdam,
Plantage Muidergracht 24,
1018 TV Amsterdam, the Netherlands\endaddress
\email koelink\@fwi.uva.nl\endemail
\address Vakgroep Toegepaste Wiskunde en Informatica,
Universiteit Gent, Krijgslaan 281-S9, B-9000 Gent,
Belgium\endaddress
\email Joris.VanderJeugt\@rug.ac.be\endemail
\date July 3, 1996\enddate
\thanks First author is supported by the Netherlands
Organization for Scientific Research (NWO)
under project number 610.06.100.
Second author is a Senior Research Associate of the National
Fund for Scientific Research of Belgium (NFWO). \endthanks
\abstract The interpretation of the Meixner-Pollaczek,
Meixner and Laguerre polynomials as overlap coefficients
in the positive discrete series representations of the
Lie algebra ${\frak{su}}(1,1)$ and the Clebsch-Gordan
decomposition leads to
generalisations of the convolution identities for
these polynomials. Using the Racah coefficients
convolution identities for continuous Hahn, Hahn and Jacobi
polynomials are obtained. From the quantised 
universal enveloping algebra
for ${\frak{su}}(1,1)$ convolution identities for
the Al-Salam and Chihara polynomials
and the Askey-Wilson polynomials are derived by using the
Clebsch-Gordan and Racah coefficients. For the quantised
universal enveloping algebra for ${\frak{su}}(2)$ $q$-Racah
polynomials are interpreted as Clebsch-Gordan coefficients,
and the linearisation coefficients for a two-parameter
family of Askey-Wilson polynomials are derived.
\endabstract
\endtopmatter
\document

\head\newsection Introduction\endhead

The representation theory of Lie algebras and quantum
algebras, or quantised universal enveloping algebras
\cite{\CharP}, is intimately linked to special functions
of (basic) hypergeometric type, see
e.g. \cite{\VileK}, \cite{\CharP}. In this paper we consider
especially the Lie algebra ${\frak{su}}(1,1)$ and its
quantum analogue $\Unc$, and we derive convolution identities
for certain orthogonal polynomials which occur as overlap
coefficients. The idea, which is due to Granovskii and Zhedanov
\cite{\GranZ}, see also \cite{\VdJeug}, is to consider
(generalised) eigenvectors of a suitable element of
the Lie algebra which is a recurrence operator
in an irreducible representation of this Lie algebra.
Then there is a relation between these eigenvectors and
the eigenvectors of this Lie algebra element
in the $n$-fold tensor product of irreducible representations of
the Lie algebra. From the tensor product decomposition in
irreducible representations for $n=2$, $3$, we obtain identities
for these eigenvectors involving Clebsch-Gordan and Racah
coefficients. In particular, if the overlap coefficients are
known in terms of special functions, we obtain identities
for these special functions in this way.

For the Lie algebra ${\frak{su}}(1,1)$ and the positive discrete
series representations a special case of this
approach is contained in Granovskii and Zhedanov
\cite{\GranZ}, but the result is
not worked out in detail. Elaborating the method of Granovskii
and Zhedanov, Van der Jeugt \cite{\VdJeug} obtains a
generalisation of the classical convolution identity
for the Laguerre polynomials
\cite{\HTFtwee, 10.12(41)}. Van der Jeugt \cite{\VdJeug} also
considers the boson Lie algebra ${\frak b}(1)$, a central
extension of the oscillator algebra, leading to a generalisation
of the convolution identity for Hermite polynomials
\cite{\HTFtwee, 10.13(38)}. The last identity follows from the
previous one by a well-known limit transition of Laguerre
polynomials to Hermite polynomials, see e.g. \cite{\KoekS}.

Apart from the Laguerre and the Hermite polynomials, also
the Meixner-Pollaczek, Meixner and Charlier polynomials,
which all fit into the Askey scheme of hypergeometric
orthogonal polynomials \cite{\AskeW}, \cite{\KoekS},
satisfy a convolution identity of the
same form. This is a straightforward consequence of the
existence of a generating function of a special kind, see
Al-Salam \cite{\AlSa}. It is also known that the Meixner-Pollaczek
and the Meixner polynomials can be interpreted as overlap
coefficients in the positive discrete series representations
of ${\frak{su}}(1,1)$, see Masson and Repka \cite{\MassR}.
In \S 3 we show how the method of Granovskii and Zhedanov
for the two-fold tensor product of positive discrete
series representations of ${\frak{su}}(1,1)$
leads to a generalisation of the
convolution identity for Meixner-Pollaczek polynomials,
{}from which generalised convolution
formulas for Meixner, Laguerre, Charlier and
Hermite polynomials can be obtained by substitution or
by limit transitions. Next using the three-fold
tensor product representation we obtain a
very general convolution identity for
continuous Hahn polynomials, and similarly for
the Hahn and Jacobi polynomials. These identities can also
be viewed as yielding connection coefficients between two sets
of orthogonal polynomials in two variables with respect
to the same orthogonality measure. With this point of view,
this result coincides with Dunkl's results \cite{\DunkPJM},
\cite{\DunkCJM}. Our derivation gives an intrinsic
explanation for the occurrence of balanced ${}_4F_3$-series
as connection coefficients; they are Racah coefficients.
Actually, the interpretation as orthogonal polynomials
in two variables works in general, and is an intrinsic
way to determine the $S$-functions in \cite{\GranZ} and
\cite{\VdJeug} in terms of orthogonal polynomials
instead of reducing a triple sum to a single sum.

In \S 4 we apply the same idea to the quantised universal
enveloping algebra $\Unc$ and its positive discrete
series representations. Due to the non-cocommutativity
of the comultiplication, which is needed to define
the tensor product representation, the tensor product
of eigenvectors is no longer an eigenvector
in the tensor product representation. This can be solved if
we restrict to operators related to so-called
twisted primitive elements in $\Unc$, see e.g.
\cite{\KoorZSE}, \cite{\KoelAAM}. Then the whole
machinery works and we obtain a generalisation
of the Al-Salam and Chihara \cite{\AlSaC} convolution
identity for the Al-Salam and Chihara polynomials by
considering the Clebsch-Gordan coefficients in the
two-fold tensor product. Going to the three-fold
tensor product representations yields a very
general convolution identity for Askey-Wilson
polynomials also involving $q$-Racah polynomials, and
Theorem~4.10 
is the key result of this paper.
Overlap coefficients are also considered in somewhat more
generality in Klimyk and Kachurik \cite{\KlimK}, but
we have to restrict ourselves to the twisted primitive
elements in order to keep the action in the tensor
product representations manageable.

It is interesting to note that in this derivation we
have a natural interpretation of the continuous Hahn,
Hahn and Jacobi polynomials as Clebsch-Gordan coefficients
for the Lie algebra ${\frak{su}}(1,1)$. Similarly, we have
an interpretation of the Askey-Wilson polynomials as
Clebsch-Gordan coefficients for the
quantised universal enveloping algebra
$\Unc$. In \S 5 we shortly discuss the corresponding
result for the quantised universal enveloping
algebra $\Usu$, where the $q$-Racah
polynomials then occur as Clebsch-Gordan coefficients.
This case can be obtained formally from the results
for $\Unc$.
Since in the dual Hopf $\ast$-algebra the so-called
zonal spherical elements are known in terms of
a two-parameter family of Askey-Wilson polynomials,
cf. \cite{\KoorZSE}, we obtain the explicit
linearisation coefficients for this subfamily of
the Askey-Wilson polynomials.

It should be remarked that there does not seem to be
an appropriate $q$-analogue of the boson Lie algebra
${\frak b}(1)$. Either, the Hopf $\ast$-algebra
structure is lacking, or, as in \cite{\KalnMM},
the recurrence in the two-fold tensor product
representation seems unmanageable.

Instead of using generalised eigenvectors we use the
spectral theory of Jacobi matrices, which we recall
briefly in \S 2. In particular we use this theory to
interpret certain recurrence operators in
$\Hi^{\otimes n}$, $n=1,2,3$, as multiplication
operators in certain weighted $L^2$-spaces on $\R^n$.
This approach exploits the theory of orthogonal
polynomials, cf. Propositions~3.3 and 4.3. 

The notation for (basic) hypergeometric series is the
standard one as in Gasper and Rahman \cite{\GaspR}.
Unexplained notions for quantised universal enveloping
algebras can be found in Chari and Pressley \cite{\CharP}.

\demo{Acknowledgement} We thank Tom Koornwinder for
useful comments. The first author thanks the Universiteit
Gent for its hospitality.
\enddemo

\head\newsection Jacobi matrices and orthogonal
polynomials\endhead

We recall some of the results on the spectral theory
of Jacobi matrices and the relation with orthogonal polynomials.
For more information we refer to Berezanski\u\i\ \cite{\Bere,
Ch.~VII, \S 1}, see also Masson and Repka \cite{\MassR},
Klimyk and Kachurik \cite{\KlimK}. The operator $J$ acting on
the standard orthonormal basis $\{ e_n\mid n\in\Zp\}$ of $\Hi$ by
$$
Je_n =a_n\, e_{n+1} + b_n\, e_n + a_{n-1}\, e_{n-1},
\qquad a_n>0,\ b_n\in\R,
\tag\eqname{\vgldefJacobi}
$$
is called a Jacobi matrix. This operator is symmetric, and its
deficiency indices are $(0,0)$ or $(1,1)$. In particular, if
the coefficients $a_n$ and $b_n$ are bounded, $J$ is
a bounded operator on $\Hi$ and thus self-adjoint. $J$ is an
unbounded self-adjoint operator
if $\sum_{n=0}^\infty a_n^{-1}=\infty$ by Carleman's condition.
Then $e_0$ is a cyclic vector for $J$, i.e. the span of finite
linear combinations of the form $J^pe_0$, $p\in\Zp$, is dense
in $\Hi$. This is the case for all Jacobi matrices considered
in this paper.

Assuming this, we can use the same coefficients $a_n$, $b_n$ to
generate polynomials $p_n(x)$ of degree $n$ in $x$ by the
recurrence relation
$$
xp_n(x) =a_n\, p_{n+1}(x) + b_n\, p_n(x) +
a_{n-1}\, p_{n-1}(x),\qquad
p_{-1}(x)=0,\ p_0(x)=1.
\tag\eqname{\vgldeforthonormpols}
$$
By Favard's theorem there exists a positive measure $m$
on the real line such that the polynomials
$p_n(x)$ are orthonormal;
$$
\int_\R p_n(x)p_m(x) \, dm(x) = \de_{n,m}.
$$
The measure is obtained by $m(B)=\langle
E(B)e_0,e_0\rangle$, $B$ Borel
set, where $E$ denotes the spectral
decomposition of the self-adjoint
operator $J$.

We can represent the operator $J$ as a multiplication
operator $M_x$ on $L^2(m)$, where $M_xf(x)=xf(x)$. For
this we define
$$
\Lambda\colon\Hi\to L^2(m), \qquad\qquad
\bigl( \Lambda e_n\bigr) (x) = p_n(x),
$$
then $\Lambda$ is a unitary operator, since
it maps an orthonormal basis onto
an orthonormal basis. Note that we use here
that the polynomials are dense in
$L^2(m)$, since the self-adjointness of $J$
implies that the corresponding
moment problem is determined. From  \thetag{\vgldefJacobi} and
\thetag{\vgldeforthonormpols} it follows that
$\Lambda \circ J = M_x \circ \Lambda$,
so that $\Lambda$ intertwines the Jacobi matrix $J$
on $\Hi$ with the multiplication
operator $M_x$ on $L^2(m)$.

\head\newsection The case ${\frak{su}}(1,1)$\endhead

The Lie algebra ${\frak{su}}(1,1)$ is given by
$$
[H,B] = 2B, \qquad [H,C] = -2C, \qquad [B,C]=H.
$$
There is a $\ast$-structure by $H^\ast=H$ and $B^\ast=-C$.

The positive discrete series representations $\pi_k$ of
${\frak{su}}(1,1)$ are unitary representations labelled by $k>0$.
The representation space is $\Hi$ equipped
with orthonormal basis $\{ e^k_n\}_{\{ n\in\Zp\} }$. The action
is given by
$$
\aligned
\pi_k(H)\, e^k_n &= 2(k+n) \, e^k_n, \\
\pi_k(B)\, e^k_n &= \sqrt{(n+1)(2k+n)}\, e^k_{n+1},\\
\pi_k(C)\, e^k_n &= -\sqrt{n(2k+n-1)}\, e^k_{n-1}.
\endaligned
\tag\eqname{\vgldefposdiscrserreps}
$$
The tensor product of two positive discrete series
representations decomposes as
$$
\pi_{k_1}\otimes \pi_{k_2} = \bigoplus_{j=0}^\infty
\pi_{k_1+k_2+j}.
\tag\eqname{\vgltensorprodddec}
$$
The corresponding intertwining operator can be expressed
by means of the Clebsch-Gordan coefficients
$$
e^k_n = \sum_{n_1,n_2} C^{k_1,k_2,k}_{n_1,n_2,n}\,
e^{k_1}_{n_1}\otimes e^{k_2}_{n_2}.
\tag\eqname{\vgldefCGC}
$$
Later we also use the notation $e^{(k_1k_2)k}_n$ for $e^k_n$ to
stress the fact that this vector arises from the decomposition
$\pi_{k_1}\otimes\pi_{k_2}$ into irreducible representations.
The Clebsch-Gordan coefficients are non-zero only if $n_1+n_2=n+j$,
$k=k_1+k_2+j$ for $j,n_1,n_2,n\in\Zp$ by considering the action of
$H$ on both sides. We normalise the Clebsch-Gordan coefficients
by $\langle e^k_0, e^{k_1}_0\otimes e^{k_2}_j\rangle >0$.

For the above results Vilenkin and Klimyk \cite{\VileK, \S 8.7}
can be consulted.

\subhead\newsubsection Clebsch-Gordan coefficients and
orthogonal polynomials \endsubhead
The Meixner-Polla\-czek polynomials are defined by
$$
P_n^{(\la)}(x;\phi) = {{(2\la)_n}\over{n!}} e^{in\phi}
\, {}_2F_1\left( {{-n,\la+ix}\atop{2\la}};1-e^{-2i\phi}\right).
\tag\eqname{\vgldefMeixnerPollaczekpols}
$$
For $\la>0$ and $0<\phi<\pi$ these are orthogonal polynomials
with respect to a positive measure on $\R$, see \cite{\KoekS}.
The orthonormal Meixner-Pollaczek polynomials
$$
p_n(x) = p_n^{(\la)}(x;\phi) = \sqrt{ {{n!}\over{\Gamma(n+2\la)}}}
P_n^{(\la)}(x;\phi)
$$
satisfy the three-term recurrence relation
$$
\align
2x\, \sin\phi\, p_n(x) &= a_n\, p_{n+1}(x) - 2(n+\la)\cos\phi\,
p_n(x) + a_{n-1}\, p_{n-1}(x),\\ a_n &= \sqrt{(n+1)(n+2\la)}.
\endalign
$$
The orthogonality measure for Meixner-Pollaczek polynomials is
absolutely continuous. Define
$$
w^{(\la)}(x;\phi) = {{ (2\sin\phi)^{2\la}}\over{2\pi}}
e^{(2\phi-\pi)x} \bigl\vert \Gamma(\la+ix)\bigr\vert^2,
$$
then
$$
\int_\R p_n^{(\la)}(x;\phi)p_m^{(\la)}(x;\phi)\,
w^{(\la)}(x;\phi)dx =\de_{nm}.
$$

Define the self-adjoint element in ${\frak{su}}(1,1)$;
$$
X_\phi = -\cos\phi\, H + B - C.
\tag\eqname{\vgldefXphi}
$$

\proclaim{Proposition \theoremname{\propspectXphi}}
$\Lambda\colon \Hi \to L^2(\R,w^{(k)}(x;\phi)dx)$,
$e^k_n\mapsto p_n^{(k)}(\cdot;\phi)$,
is a unitary mapping
intertwining $\pi_k(X_\phi)$ acting in $\Hi$ with $M_{2x\sin\phi}$
on $L^2(\R,w^{(k)}(x;\phi)dx)$.
\endproclaim

Here, and elsewhere, $M_g$ denotes multiplication by the function
$g$, so $M_gf(x)=g(x)f(x)$.

\demo{Proof} Use \thetag{\vgldefposdiscrserreps} and
\thetag{\vgldefXphi} to see that $\pi_k(X_\phi)$ is a Jacobi
matrix. Next compare the coefficients with the three-term
recurrence relation for the orthonormal Meixner-Pollaczek
polynomials to find the result as in \S 2. 
\qed\enddemo

Proposition \thmref{\propspectXphi}
states that $v^k(x)=\sum_{n=0}^\infty p_n^{(k)}(x;\phi)\, e^k_n$
is a generalised eigenvector for $\pi_k(X_\phi)$ for the
eigenvalue $2x\sin\phi$. Next we study the action of $X_\phi$
in the tensor product representation $\pi_{k_1}\otimes\pi_{k_2}$.
Recall that $\De(X_\phi)=1\otimes X_\phi+X_\phi\otimes 1$.

\proclaim{Proposition \theoremname{\propspectDeltaXphi}}
$\Upsilon\colon \Hi\otimes\Hi \to
L^2(\R^2,w^{(k_1)}(x_1;\phi)w^{(k_2)}(x_2;\phi)dx_1dx_2)$,
defined by $e^{k_1}_{n_1}\otimes e^{k_2}_{n_2} \mapsto
p_{n_1}^{(k_1)}(x_1;\phi)p_{n_2}^{(k_2)}(x_2;\phi)$
is a unitary mapping intertwining
$\pi_{k_1}\otimes \pi_{k_2}(\De(X_\phi))$
with $M_{2(x_1+x_2)\sin\phi}$.
\endproclaim

\demo{Proof} This can be seen by using the mapping $\Lambda$ of
Proposition \thmref{\propspectXphi} in the second tensor factor
and solving the resulting three-term recurrence in the first
factor.
\qed\enddemo

Proposition \thmref{\propspectDeltaXphi} states that
$$
v^{k_1,k_2}(x_1,x_2) = \sum_{n_1,n_2=0}^\infty
p_{n_1}^{(k_1)}(x_1;\phi) p_{n_2}^{(k_2)}(x_2;\phi)\
e^{k_1}_{n_1}\otimes e^{k_2}_{n_2}
$$
are generalised eigenvectors for
$\pi_{k_1}\otimes\pi_{k_2}(\De(X_\phi))$
for the eigenvalue $2(x_1+x_2)\sin\phi$.

So $\Upsilon$ maps the basis $e^{k_1}_{n_1}\otimes e^{k_2}_{n_2}$
onto orthonormal polynomials in two variables. By the
Clebsch-Gordan decomposition \thetag{\vgltensorprodddec}
there exists another orthonormal basis $e_n^k$ for the
tensor product representation space. So $\Upsilon e^k_n$
gives another set of orthonormal polynomials in two variables
in $L^2(\R^2,w^{(k_1)}(x_1;\phi)w^{(k_2)}(x_2;\phi)dx_1dx_2)$.

In order to formulate the
result we need the continuous Hahn polynomials
\cite{\KoekS}, \cite{\KoelPAMS};
$$
p_n(x;a,b,c,d)= i^n {{(a+c)_n(a+d)_n}\over{n!}}\,
{}_3F_2\left( {{-n,n+a+b+c+d-1,a+ix}\atop{a+c,\ a+d}};1\right)
\tag\eqname{\vgldefcontHahnpols}
$$
satisfying the orthogonality relations for $\Re (a,b,c,d)>0$
$$
\multline
{1\over{2\pi}} \int_\R \Gamma(a+ix)\Gamma(b+ix)
\Gamma(c-ix)\Gamma(d-ix)\, p_n(x;a,b,c,d)p_m(x;a,b,c,d)\, dx = \\
\de_{nm} {{\Gamma(n+a+c)\Gamma(n+a+d)\Gamma(n+b+c)\Gamma(n+b+d)}
\over {n!\, (2n+a+b+c+d-1) \Gamma(n+a+b+c+d-1)}}.
\endmultline
$$
The orthogonality measure is positive for $a=\bar c$, $b=\bar d$.

\proclaim{Proposition \theoremname{\propspectdiffbasisDeltaXphi}}
In $L^2(\R^2,w^{(k_1)}(x_1;\phi)w^{(k_2)}(x_2;\phi)dx_1dx_2)$
we have
$$
\align
\Upsilon e^k_n (x_1,x_2) &= p_n^{(k)}(x_1+x_2;\phi) \, \Upsilon
e^k_0 (x_1,x_2), \\
\Upsilon e^k_0 (x_1,x_2) &= C p_j(x_1;k_1,
k_2-i(x_1+x_2),k_1,k_2+i(x_1+x_2)),\\
C &= (-2\sin\phi)^j
\sqrt{ {{j!\, (2j+2k_1+2k_2-1)\Gamma(j+2k_1+2k_2-1)}\over
{\Gamma(2k_1+j) \Gamma(2k_2+j)}} }.
\endalign
$$
\endproclaim

Note that $\Upsilon e^k_0 (x_1,x_2)$ is indeed a polynomial in
$x_1$, $x_2$.

\demo{Proof} The first statement follows from use of the
intertwining of Proposition \thmref{\propspectDeltaXphi} and the
intertwining of \thetag{\vgltensorprodddec};
$$
2(x_1+x_2)\sin\phi\,  \Upsilon e^k_n(x_1,x_2) =
M_{2(x_1+x_2)\sin\phi} \Upsilon e^k_n(x_1,x_2) =
\bigl( \Upsilon \pi_k(X_\phi)\, e^k_n\bigr) (x_1,x_2),
$$
which gives a three-term recurrence relation for $\Upsilon e^k_n$
with respect to $n$ of the same form as in Proposition
\thmref{\propspectXphi}. Taking into account the initial conditions
proves the first statement.

The prove the second statement we note that for $k=k_1+k_2+j$,
$l=k_1+k_2+i$,
$$
\multline
\de_{ij}\de_{mn} = \langle e^k_n,e^l_m\rangle =
\langle \Upsilon e^k_n, \Upsilon e^l_m\rangle =
\iint_{\R^2} p_n^{(k)}(x_1+x_2;\phi)
p_m^{(l)}(x_1+x_2;\phi) \\ \times \Bigl( \Upsilon e^k_0 (x_1,x_2)
\Upsilon e^l_0 (x_1,x_2)\Bigr)
w^{(k_1)}(x_1;\phi)w^{(k_2)}(x_2;\phi)\, dx_1dx_2,
\endmultline
$$
by the first statement and Proposition \thmref{\propspectDeltaXphi}.
Introduce $s=x_1+x_2$, $t=x_1$, then we find
$$
\de_{ij}\de_{mn} = \int_\R p_n^{(k)}(s;\phi) p_m^{(l)}(s;\phi)
\int_\R \Upsilon e^k_0 (t,s-t)
\Upsilon e^l_0 (t,s-t)
w^{(k_1)}(t;\phi)w^{(k_2)}(s-t;\phi)\, dt ds
$$
In case $k=l$, or $i=j$, we see that the inner integral must equal
the normalised orthogonality measure for the Meixner-Pollaczek
polynomials $p_n^{(k)}(s;\phi)$, since the corresponding moment
problem is determined. In case $k\not=l$, or $i\not=j$, we
conclude that the inner integral integrated against any polynomial
gives zero, so that it must be zero since the polynomials are
dense in
$L^2(\R^2,w^{(k_1)}(x_1;\phi)w^{(k_2)}(x_2;\phi)dx_1dx_2)$.
So we get
$$
\multline
\de_{ij} w^{(k)}(s;\phi) =
e^{(2\phi-\pi)s} {{(2\sin\phi)^{2k_1+2k_2}}\over{4\pi^2}}
\int_\R \Upsilon e^k_0 (t,s-t)
\Upsilon e^l_0 (t,s-t) \\
\times
\Gamma(k_1+it)\Gamma(k_2-is+it)
\Gamma(k_1-it)\Gamma(k_2+is-it) \, dt.
\endmultline
$$
Apply $\Upsilon$ to \thetag{\vgldefCGC} for $n=0$ to see that
$\Upsilon e^k_0 (t,s-t)$ is a polynomial of degree $j$ in $t$.
Hence, $\Upsilon e^k_0 (t,s-t)$ is a multiple of a continuous
Hahn polynomial of degree $j$ with the parameters as in the
proposition.

The value of the constant follows from comparing the squared
norms up to a sign. The sign is determined from the condition
on the Clebsch-Gordan coefficients. This implies
$0<\langle \Upsilon e^k_0,
\Upsilon e^{k_1}_0\otimes e^{k_2}_j\rangle$
and using the first two parts of the proposition and Proposition
\thmref{\propspectDeltaXphi} shows that the sign of $C$ follows
{}from the sign of a double integral of two orthogonal polynomials.
Only the integral over $x_2$ is relevant, and the sign of $C$
equals the sign of the leading coefficient of the continuous Hahn
polynomials viewed as a polynomial in $x_2$, which is $(-1)^j$.
\qed\enddemo

So we can now apply $\Upsilon$ to \thetag{\vgldefCGC} to find,
$k=k_1+k_2+j$,
$$
\multline
\sum_{n_1+n_2=n+j} C^{k_1,k_2,k}_{n_1,n_2,n}\,
p_{n_1}^{(k_1)}(x_1;\phi) p_{n_2}^{(k_2)}(x_2;\phi) =
C p_n^{(k)}(x_1+x_2;\phi) \\ \times
p_j(x_1;k_1,
k_2-i(x_1+x_2),k_1,k_2+i(x_1+x_2)).
\endmultline
\tag\eqname{\vgladdformoneMP}
$$
The Clebsch-Gordan coefficients remain to be determined, and this
can be done from this formula, see \cite{\VdJeug}. They can be
expressed in terms of ${}_3F_2$-series, which are known as
Hahn polynomials. Using the Hahn polynomials defined by
$$
Q_n(x;a,b,N) = {}_3F_2\left(
{{-n,n+a+b+1,-x}\atop{a+1,\ -N}};1\right)
\tag\eqname{\vgldefHahnpols}
$$
for $N\in\Zp$, $0\leq n\leq N$, we have, with $k=k_1+k_2+j$,
$n_1+n_2=n+j$,
$$
\multline
C^{k_1,k_2,k}_{n_1,n_2,n} = \sqrt{
{{(2k_1)_{n_1}(2k_2)_{n_2}(2k_1)_j}\over
{n!n_1!n_2!j!\, (2k_1+2k_2+2j)_n (2k_2)_j (2k_1+2k_2+j-1)_j}} }\\
\times (n+j)!\, Q_j(n_1;2k_1-1,2k_2-1;n+j),
\endmultline
$$
see \cite{\VileK, \S 8.7} for another proof.

Using this in \thetag{\vgladdformoneMP} gives an identity
in a weighted $L^2$-space, but since it
is a polynomial identity it holds for all $x_1$, $x_2$.
Simplifying proves the following theorem.

\proclaim{Theorem \theoremname{\thmgenconvoidforMP}} With the
notation for continuous Hahn, Meixner-Pollaczek and Hahn
polynomials as in \thetag{\vgldefMeixnerPollaczekpols},
\thetag{\vgldefcontHahnpols} and \thetag{\vgldefHahnpols}
the following convolution formula holds:
$$
\multline
{{n+j}\choose n}
\sum_{l=0}^{n+j} Q_j(l;2k_1-1,2k_2-1,n+j)\, P_l^{(k_1)}(x_1;\phi)
\,  P_{n+j-l}^{(k_2)}(x_2;\phi) = \\
{{(-2\sin\phi)^j}\over{(2k_1)_j}}
\, P_n^{(k_1+k_2+j)}(x_1+x_2;\phi)\,
 p_j(x_1;k_1,
k_2-i(x_1+x_2),k_1,k_2+i(x_1+x_2)).
\endmultline
$$
\endproclaim

\demo{Remark \theoremname{\remthmgenconvoidforMP}} (i) The case
$j=0$ gives back the convolution identity for the
Meixner-Pollaczek polynomials, see e.g. \cite{\AlSa, \S 8},
\cite{\AlSaC}. The case $n=0$ gives another convolution
identity for Meixner-Pollaczek polynomials, since the Hahn
polynomial reduces to a summable ${}_2F_1$-series.

(ii) Note that the polynomials on both sides of the formula in
Theorem \thmref{\thmgenconvoidforMP} are orthogonal polynomials
in two variables for the space
$L^2(\R^2,w^{(k_1)}(x_1;\phi)w^{(k_2)}(x_2;\phi)dx_1dx_2)$,
so we have proved a connection coefficient formula for these
polynomials. The dual connection coefficient formula follows from
the orthogonality of the Clebsch-Gordan matrix, or equivalently,
{}from the orthogonality relations for the dual Hahn polynomials.

(iii) Theorem \thmref{\thmgenconvoidforMP} shows that the
continuous Hahn polynomials have an interpretation as
Clebsch-Gordan coefficients for ${\frak{su}}(1,1)$. Using the
generalised eigenvectors we formally have, cf. \thetag{\vgldefCGC},
$$
v^{k_1,k_2}(x_1,x_2) = \sum_k C\,p_j(x_1;k_1,
k_2-i(x_1+x_2),k_1,k_2+i(x_1+x_2))\, v^k(x_1+x_2),
$$
with $C$ as in Proposition \thmref{\propspectdiffbasisDeltaXphi}.
The dual relations can be written using the orthogonality measure
for the continuous Hahn polynomials.
\enddemo

Recall the definition of the Laguerre polynomials
$L_n^{(a)}(x) = (a+1)_n/n!\, {}_1F_1(-n;a+1;x)$,
the Jacobi polynomials
$P_n^{(a,b)}(x) = (a+1)_n/n!\, {}_2F_1(-n,n+a+b+1;a+1;(1-x)/2)$
and the Meixner polynomials $M_n(x;\be;c) =
{}_2F_1(-n,-x;\be;1-c^{-1})$.

\proclaim{Corollary \theoremname{\corthmgenconvoidforMP}}
{\rm (i) (\cite{\VdJeug})}
The Laguerre polynomials satisfy the following convolution identity
$$
\multline
\sum_{l=0}^{n+j} Q_j(l;a,b,n+j)\, L_l^{(a)}(x_1)
\,  L_{n+j-l}^{(b)}(x_2) = {{(-1)^j n!\, j!}\over{(a+1)_j (n+j)!}}\\
\times L_n^{(a+b+1+2j)}(x_1+x_2) \, (x_1+x_2)^j P_j^{(a,b)}\left(
{{x_2-x_1}\over{x_1+x_2}}\right).
\endmultline
$$
{\rm (ii)}
The Meixner polynomials satisfy the following convolution identity
$$
\multline
(c^{-1}-1)^{-j} \sum_{l=0}^{n+j} {{(a)_l (b)_{n+j-l}}\over
{l!\,(n+j-l)!}} Q_j(l;a-1,b-1,n+j)\, M_l(x_1;a;c)
\,  M_{n+j-l}(x_2;b;c) = \\
 {{(a+b+2j)_n}\over{(n+j)!}}\, M_n(x_1+x_2-j;a+b+2j;c) \,
(-x_1-x_2)_j Q_j(x_1;a-1,b-1,x_1+x_2).
\endmultline
$$
\endproclaim

\demo{Proof} The first case follows from the limit transition of
the Meixner-Pollaczek polynomials to the Laguerre polynomials;
$\lim_{\phi\downarrow 0}
P_n^{((a+1)/2)}(-2x/\phi;\phi)=L_n^{(a)}(x)$.
In this limit transition the continuous Hahn polynomials tend to
the Jacobi polynomials.

The second case follows from the substitution
$\phi = \ln c/2i$, and replacing $x_1$ and $x_2$ by $ik_1+ix_1$
and $ik_2+ix_2$. For this substitution the continuous Hahn
polynomials go over into the Hahn polynomials.
\qed\enddemo

\demo{Remark \theoremname{\remcorthmgenconvoidforMP}}
(i) The case $j=0$ in both formulas
gives back the convolution identities for the
Laguerre and Meixner polynomials,
see e.g. \cite{\AlSa}, \cite{\AlSaC}, \cite{\HTFtwee, 10.12(41)},
and the case $n=0$ gives
another convolution identity for the Laguerre and Meixner
polynomials. Again these
formulas can be viewed as connection coefficient formulas
for orthogonal polynomials in two variables.

(ii) The identities of Corollary \thmref{\corthmgenconvoidforMP}
can be obtained by considering the action of $X=-H+B-C$
in the representations $\pi_k$ and $\pi_{k_1}\otimes \pi_{k_2}$
for the Laguerre case, see \cite{\VdJeug}, and by considering
the action of $X_c=-(1+c)/2\sqrt{c} H+B-C$, $0<c<1$,
in the representations $\pi_k$ and $\pi_{k_1}\otimes \pi_{k_2}$
for the Meixner case. The limit case $c\uparrow 1$ in the Meixner
result gives the  Laguerre result. In this case we can interpret
the Jacobi and Hahn polynomials as Clebsch-Gordan coefficients,
cf. Remark \thmref{\remthmgenconvoidforMP}(iii).

(iii) Corollary \thmref{\corthmgenconvoidforMP}(ii) is equivalent
to Theorem \thmref{\thmgenconvoidforMP} by the same substitution.
Theorem \thmref{\thmgenconvoidforMP} can also be obtained
{}from Corollary \thmref{\corthmgenconvoidforMP}(i) by a double
application of the Mellin transform. For this we have to use that
the Laguerre polynomials are mapped onto Meixner-Pollaczek
polynomials, cf. \cite{\KoorJMP, \S 3}, and that the Jacobi
polynomials are mapped onto the continuous Hahn polynomials,
cf. \cite{\KoelPAMS, (3.4) with $\Gamma(\be-i\la)$
replaced by $\Gamma(\be+i\la)$}.
\enddemo

The other hypergeometric orthogonal polynomials satisfying a
convolution identity are the Charlier and Hermite polynomials,
cf. \cite{\AlSa}, \cite{\AlSaC}. These identities can be obtained
by taking the appropriate limits from the Meixner polynomials
to the Charlier polynomials and from the Laguerre polynomials
to the Hermite polynomials, cf. e.g. \cite{\KoekS}. The Hahn
polynomials tend to Krawtchouk polynomials and the Jacobi
polynomials tend to Hermite polynomials. We use the notation
$K_n(x;p,N) = {}_2F_1(-n,-x;-N;p^{-1})$ for Krawtchouk
polynomials, $C_n(x;a)= {}_2F_0(-n,-x;-;a^{-1})$ for
Charlier polynomials and
$H_n(x)= (2x)^n{}_2F_0(-n/2,-(n-1)/2;-;-x^{-2})$ for the
Hermite polynomials.

\proclaim{Corollary \theoremname{\corcorthmgenconvoidforMP}}
{\rm (i) (\cite{\VdJeug})}
The Hermite polynomials satisfy the following convolution identity
$$
\multline
\sum_{l=0}^{n+j} K_j(l;{{a^2}\over{a^2+b^2}},n+j)
{{a^l}\over{l!}}H_l(x){{b^{n+j-l}}\over{(n+j-l)!}}H_{n+j-l}(y) = \\
{{(a^2+b^2)^{(n+j)/2}}\over{(n+j)!}} \left( {b\over a}\right)^j
H_n \left( {{ax+by}\over{\sqrt{a^2+b^2}}}\right)
H_j \left( {{ay-bx}\over{\sqrt{a^2+b^2}}}\right).
\endmultline
$$
{\rm (ii)} The Charlier polynomials satisfy the following
convolution identity
$$
\multline
\sum_{l=0}^{n+j}  {{n+j}\choose l} \al^l \be^{n+j-l}
K_j(l; {{\al}\over{\al+\be}},n+j) C_l(x;\al)\,
C_{n+j-l}(y;\be) = \\
(-1)^j (\al+\be)^n C_n(x+y-j;\al+\be) \, (-x-y)_j
 K_j(x;{{\al}\over{\al+\be}},x+y).
\endmultline
$$
\endproclaim

\demo{Remark \theoremname{\remcorcorthmgenconvoidforMP}} (i)
Again the case $j=0$ gives known convolution formulas,
cf. \cite{\AlSa, \S 8}, \cite{\AlSaC}, \cite{\HTFtwee, 10.13(40)}.
Corollary \thmref{\corcorthmgenconvoidforMP}(ii) is derived in a
different way in Vilenkin and Klimyk \cite{\VileK, \S 8.6.5}.

(ii) This time the identities have a similar interpretation,
but now we have to use the Lie algebra ${\frak b}(1)$, a central
extension of the oscillator algebra, cf. \cite{\VdJeug}.
In particular we can now interpret the Hermite and Charlier
polynomials as Clebsch-Gordan coefficients.
\enddemo

\subhead\newsubsection Racah coefficients and orthogonal
polynomials \endsubhead
In the tensor product of three positive discrete
series representations $\pi_{k_1} \otimes
\pi_{k_2} \otimes \pi_{k_3}$ of ${\frak{su}}(1,1)$
we consider the following orthogonal bases;
$$
\align
e^{((k_1k_2)k_{12}k_3)k}_n &= \sum_{n_{12},n_3}
C^{k_{12},k_3,k}_{n_{12},n_3,n}\ e^{(k_1k_2)k_{12}}_{n_{12}}\otimes
e^{k_3}_{n_3} \tag\eqname{\vgldeforthobasisdrieeen}\\
&= \sum_{n_1,n_2,n_3,n_{12}}C^{k_1,k_2,k_{12}}_{n_1,n_2,n_{12}}
C^{k_{12},k_3,k}_{n_{12},n_3,n} \ e^{k_1}_{n_1}
\otimes e^{k_2}_{n_2} \otimes
e^{k_3}_{n_3}, \tag\eqname{\vgldeforthobasisdrietwee}
\endalign
$$
and
$$
\align
e^{(k_1(k_2k_3)k_{23})k}_n &= \sum_{n_1,n_{23}}
C^{k_{1},k_{23},k}_{n_{1},n_{23},n} \ e^{k_1}_{n_1} \otimes
e^{(k_2k_3)k_{23}}_{n_{23}}
\tag\eqname{\vgldeforthobasisdrievier}\\
&= \sum_{n_1,n_2,n_3,n_{23}}C^{k_2,k_3,k_{23}}_{n_2,n_3,n_{23}}
C^{k_{1},k_{23},k}_{n_{1},n_{23},n} \ e^{k_1}_{n_1} \otimes
e^{k_2}_{n_2} \otimes
e^{k_3}_{n_3}. \tag\eqname{\vgldeforthobasisdrievier}
\endalign
$$
Here we use the extended notation $e^{(k_1k_2)k}_n$ for the basis
of the tensor product decomposition to keep track of how the
decomposition is obtained.

These bases are connected by the Racah coefficients, which
leads to an intertwiner for the action of ${\frak{su}}(1,1)$.
The Racah coefficients are defined by
$$
e^{((k_1k_2)k_{12}k_3)k}_n = \sum_{k_{23}}
U^{k_1,k_2,k_{12}}_{k_3,k,k_{23}} \ e^{(k_1(k_2k_3)k_{23})k}_n .
\tag\eqname{\vgldefracahcoeff}
$$
In the previous formulas the following constraints hold
$$
\gathered
k_{12}=k_1+k_2+j_{12}, \qquad k_{23}=k_2+k_3+j_{23}, \\
k=k_{12}+k_3+j =k_1+k_{23}+j',\qquad
j_{12},j,j_{23},j' \in \Zp,\ \text{and}\ j_{12}+j=j_{23}+j'.
\endgathered
\tag\eqname{\vglconventiosumvariables}
$$
Thus all above sums are finite sums.

Recall that $(1 \otimes \Delta)(\Delta(X_\phi))=
1\otimes1\otimes X_\phi +1\otimes X_\phi\otimes 1+
X_\phi\otimes 1\otimes 1$. The following proposition is proved as
Proposition \thmref{\propspectDeltaXphi}.

\proclaim{Proposition \theoremname{\propspceDeltaDeltaXphi}}
Define the unitary mapping
$$
\Theta \colon \ell^2(\Zp)\otimes \ell^2(\Zp) \otimes \ell^2(\Zp)
\rightarrow L^2(\R^3,w^{(k_1)}(x_1;\phi)
w^{(k_2)}(x_2;\phi) w^{(k_3)}(x_3;\phi)dx_1dx_2dx_3)
$$
by
$$
\Theta \colon e^{k_1}_{n_1} \otimes e^{k_2}_{n_2} \otimes
e^{k_3}_{n_3} \mapsto  p^{(k_1)}_{n_1}(x_1;\phi)
p^{(k_2)}_{n_2}(x_2;\phi) p^{(k_3)}_{n_3}(x_3;\phi),
$$
then $\Theta$ intertwines $\pi_{k_1} \otimes
\pi_{k_2} \otimes \pi_{k_3}((1 \otimes \Delta)(\Delta(X_\phi))$
with $M_{2(x_1+x_2+x_3)\sin\phi}$.
\endproclaim

\demo{Remark \theoremname{\rempropspceDeltaDeltaXphi}} Let
$\Lambda^{(k)}=\Lambda$ be the unitary mapping defined in
Proposition \thmref{\propspectXphi} and
$\Upsilon^{(k_1k_2)}=\Upsilon$
be the unitary mapping defined in Proposition
\thmref{\propspectDeltaXphi}. Using the identifications
$$
\align
&L^2(\R^3,w^{(k_1)}(x_1;\phi)
w^{(k_2)}(x_2;\phi) w^{(k_3)}(x_3;\phi)dx_1dx_2dx_3) \\
&\qquad = L^2(\R,w^{(k_1)}(x_1;\phi)dx_1) \otimes
L^2(\R^2, w^{(k_2)}(x_2;\phi) w^{(k_3)}(x_3;\phi)dx_2dx_3) \\
&\qquad = L^2(\R^2,w^{(k_1)}(x_1;\phi)
w^{(k_2)}(x_2;\phi) )dx_1dx_2) \otimes
L^2(\R, w^{(k_3)}(x_3;\phi)dx_3),
\endalign
$$
we have $\Theta= \Lambda^{(k_1)}\otimes \Upsilon^{(k_2k_3)}
=\Upsilon^{(k_1k_2)}\otimes \Lambda^{(k_3)}$. Hence, for the
orthogonal bases on the right hand side of
\thetag{\vgldeforthobasisdrieeen} and
\thetag{\vgldeforthobasisdrievier} we have
$$
\align
\Theta e^{(k_1k_2)k_{12}}_{n_{12}}\otimes
e^{k_3}_{n_3} &=
\Bigl(\Upsilon^{(k_1k_2)} e_{n_{12}}^{(k_1k_2)k_{12}}\Bigr)\,
\Bigl( \Lambda^{(k_3)} e^{k_3}_{n_3}\Bigr), \\
\Theta e^{k_1}_{n_1} \otimes
e^{(k_2k_3)k_{23}}_{n_{23}} &= \Bigl(
\Lambda^{(k_1)}e^{k_1}_{n_1}\Bigr)\,\Bigl(
\Upsilon^{(k_2k_3)} e_{n_{23}}^{(k_2k_3)k_{23}}\Bigr).
\endalign
$$
And the right hand sides are known from Propositions
\thmref\thmref{\propspectXphi} and \thmref{\propspectDeltaXphi}
in terms of Meixner-Pollaczek polynomials times continuous
Hahn polynomials.
\enddemo

\proclaim{Proposition
\theoremname{\propspectdiffbasisDeltaDeltaXphi}}
{\rm (i)} The following expressions hold;
$$
\align
\Theta (e^{((k_1k_2)k_{12}k_3)k}_n)(x_1,x_2,x_3) &=
p_n^{(k)}(x_1+x_2+x_3;\phi)\
\Theta (e^{((k_1k_2)k_{12}k_3)k}_0)(x_1,x_2,x_3), \\
\Theta (e^{((k_1k_2)k_{12}k_3)k}_0)(x_1,x_2,x_3) &=
\Upsilon^{(k_1k_2)} e^{(k_1k_2)k_{12}}_0(x_1,x_2)\
\Upsilon^{(k_{12}k_3)} e^{(k_{12}k_3)k}_0(x_1+x_2,x_3).
\endalign
$$
{\rm (ii)} The following expressions hold;
$$
\align
\Theta (e^{(k_1(k_2k_3)k_{23})k}_n) (x_1,x_2,x_3) &=
p_n^{(k)}(x_1+x_2+x_3;\phi)
\ \Theta (e^{(k_1(k_2k_3)k_{23})k}_0) (x_1,x_2,x_3), \\
\Theta (e^{(k_1(k_2k_3)k_{23})k}_0) (x_1,x_2,x_3) &=
\Upsilon^{(k_2k_3)}
e^{(k_2k_3)k_{23}}_0(x_2,x_3)\  \Upsilon^{(k_1k_{23})}
e^{(k_1k_{23})k}_0(x_1,x_2+x_3) .
\endalign
$$
\endproclaim

\demo{Proof} Statement (ii) is proved analogously as statement (i).
The first statement of (i) follows from Proposition
\thmref{\propspceDeltaDeltaXphi} and the decomposition of the
three-fold tensor product, cf. Proposition
\thmref{\propspectdiffbasisDeltaXphi}.

For the second statement we use \thetag{\vgldeforthobasisdrieeen},
Remark \thmref{\rempropspceDeltaDeltaXphi} and
Propositions \thmref{\propspectdiffbasisDeltaXphi}
and \thmref{\propspectXphi} to find
$$
\multline
\Theta (e^{((k_1k_2)k_{12}k_3)k}_0)(x_1,x_2,x_3)
=\\ \Bigl(\Upsilon^{(k_1k_2)}e^{(k_1k_2)k_{12}}_0\Bigr)(x_1,x_2)
\sum_{n_{12}+n_3=j}
C^{k_{12},k_3,k}_{n_{12},n_3,0}\
p^{(k_{12})}_{n_{12}}(x_1+x_2;\phi)
p^{(k_3)}_{n_3}(x_3;\phi).
\endmultline
$$
The sum can be evaluated as $\Bigl(\Upsilon^{(k_{12}k_3)}
e^{(k_{12}k_3)k}_0\Bigr)(x_1+x_2,x_3)$ by
\thetag{\vgladdformoneMP}.
\qed\enddemo

Next we apply $\Theta$ to \thetag{\vgldefracahcoeff}, then it
follows from Proposition
\thmref{\propspectdiffbasisDeltaDeltaXphi} that
we can divide both sides by the Meixner-Pollaczek polynomial
of degree $n$. Since $\Theta$ is unitary we obtain
the Wigner-Eckart theorem, stating that the
Racah coefficients in \thetag{\vgldefracahcoeff} are
independent of $n$. So we can restrict
to the case $n=0$ of \thetag{\vgldefracahcoeff} before applying
$\Theta$ without loss of generality. We obtain
$$
\multline
\sum_{j_{23}}
U^{k_1,k_2,k_{12}}_{k_3,k,k_{23}} \, \Bigl(\Upsilon^{(k_2k_3)}
e^{(k_2k_3)k_{23}}_0\Bigr) (x_2,x_3) \, \Bigl(
\Upsilon^{(k_1k_{23})} e^{(k_1k_{23})k}_0\Bigr)(x_1,x_2+x_3) =\\
\Bigl(\Upsilon^{(k_1k_2)} e^{(k_1k_2)k_{12}}_0\Bigr)(x_1,x_2)
\, \Bigl(\Upsilon^{(k_{12}k_3)}
e^{(k_{12}k_3)k}_0\Bigr)(x_1+x_2,x_3).
\endmultline
\tag\eqname{\vglgeneralidRCone}
$$

The Racah coefficients remain to be determined, and this can
actually be done from \thetag{\vglgeneralidRCone}, see
\cite{\VdJeug}. One can either copy the expression
\cite{\VdJeug, (4.8)}, or use the limit $q\uparrow 1$ of the
expression for the $q$-Racah coefficient given in
Proposition~4.9. 
Both lead to the following expression of
the Racah coefficients in terms of balanced ${}_4F_3$-series;
$$
\multline
U^{k_1,k_2,k_{12}}_{k_3,k,k_{23}}=
\left(j+j_{12}\atop j_{23}\right)
{(2k_2)_{j_{12}}(2k_3)_j(2k_1+2k_2+2k_3+j+j_{12}-1)_{j_{23}}
\over (2k_3,2k_2+2k_3+j_{23}-1)_{j_{23}}
(2k_2+2k_3+2j_{23})_{j'} } \\
\times \left({j'!(2k_1,2k_{23},2k_1+2k_{23}+j'-1)_{j'} \;
j_{23}!(2k_2,2k_3,2k_2+2k_3+j_{23}-1)_{j_{23}} \over
j!(2k_{12},2k_3,2k_{12}+2k_3+j-1)_j\;
j_{12}!(2k_1,2k_2,2k_1+2k_2+j_{12}-1)_{j_{12}} }
\right)^{1/2} \\
\times {}_4F_3\left(
{2k_1+2k_2+j_{12}-1,2k_2+2k_3+j_{23}-1,-j_{12},-j_{23} \atop
2k_2,2k_1+2k_2+2k_3+j+j_{12}-1,-j-j_{12}}; 1\right),
\endmultline
\tag\eqname{\vglRacahCasbalancedhypser}
$$
with the convention \thetag{\vglconventiosumvariables}.

The Racah coefficients can be rewritten in terms of the
Racah polynomials defined by
$$
R_n(\la(x);\al,\be,\ga,\de)= {}_4F_3\left(
{{-n, n+\al+\be+1,-x,x+\ga+\de+1}\atop
{\al+1,\ \be+\de+1,\ \ga+1}};1 \right),
\tag\eqname{\vgldefRacahpols}
$$
where $\la(x)=x(x+\ga+\de+1)$, one of the lower parameters
equals $-N$, $N\in\Zp$ and $0\leq n\leq N$, cf. \cite{\KoekS}.
The orthogonality
relations for the Racah polynomials follow from the fact that the
Racah coefficients form a unitary matrix.

So we obtain the following theorem by simplifying
\thetag{\vglgeneralidRCone} using
$s=x_1+x_2+x_3$ and the explicit expression
\thetag{\vglRacahCasbalancedhypser}.

\proclaim{Theorem \theoremname{\thmgeneralconvoforcontinuousHahn}}
The continuous Hahn polynomials satisfy the following convolution
identity
$$
\align
\sum_{l=0}^{n+j}&
\left(j+n\atop n\right)
{(2k_2)_n(2k_3)_j(2k_1+2k_2+2k_3+j+n-1)_l
\over (2k_3)_l(2k_2+2k_3+l-1)_l (2k_2+2k_3+2l)_{j+n-l} }\\
&\times R_l(\la(n);2k_2-1,2k_3-1,-j-n-1,2k_1+2k_2+j+n-1) \\
&\times
p_{n+j-l}(x_1;k_1,k_2+k_3+l-is,k_1,k_2+k_3+l+is)\\ &\times
p_l(x_2;k_2,k_3-i(s-x_1),k_2,k_3+i(s-x_1))\\
=&\, p_n(x_1;k_1,k_2-i(x_1+x_2),k_1,k_2+i(x_1+x_2))
\\ &\times
p_j(x_1+x_2;k_1+k_2+n,k_3-is,k_1+k_2+n,k_3+is),
\endalign
$$
with the notation as in \thetag{\vgldefcontHahnpols},
\thetag{\vgldefRacahpols}.
\endproclaim

\demo{Remark \theoremname{\remthmgeneralconvoforcontinuousHahn}}
(i) Theorem \thmref{\thmgeneralconvoforcontinuousHahn} can be
considered as a connection coefficient formula between two systems
of orthogonal polynomials for the orthogonality measure
$$
\Gamma(k_1+ix_1)\Gamma(k_1-ix_1)\Gamma(k_2+ix_2)\Gamma(k_2-ix_2)
\Gamma(k_3+i(s-x_1-x_2))\Gamma(k_3-i(s-x_1-x_2))\, dx_1dx_2
$$
on $\R^2$. This follows from substituting $s=x_1+x_2+x_3$
in the weighted $L^2$-space of Proposition
\thmref{\propspceDeltaDeltaXphi} and leaving out
the integration with respect to $s$, which can be done
by Proposition \thmref{\propspectdiffbasisDeltaDeltaXphi}
and the Wigner-Eckart theorem.

(ii) Theorem \thmref{\thmgenconvoidforMP} can be obtained as
a limit case of Theorem \thmref{\thmgeneralconvoforcontinuousHahn}
by letting $k_3\to\infty$ and using the limit transition
of the continuous Hahn polynomials to the Meixner-Pollaczek
polynomials, see e.g. \cite{\KoekS}. Note that
Theorem \thmref{\thmgenconvoidforMP} is used in the
derivation of Theorem \thmref{\thmgeneralconvoforcontinuousHahn}.

(iii) Application of $\Theta$ to
\thetag{\vgldeforthobasisdrieeen}--
\thetag{\vgldeforthobasisdrievier} gives results which are
immediately derivable from Theorem \thmref{\thmgenconvoidforMP}.
\enddemo

\proclaim{Corollary
\theoremname{\corthmgeneralconvoforcontinuousHahn}}
{\rm (i) (\cite{\VdJeug})} The Jacobi polynomials satisfy the
convolution identity
$$
\align
\sum_{l=0}^{n+j}&
\left(j+n\atop n\right)
{(b+1)_n(c+1)_j(a+b+c+j+n+2)_l
\over (c+1)_l(b+c+l+1)_l (b+c+2l+2)_{j+n-l} }\\
&\times R_l(\la(n);b,c,-j-n-1,a+b+j+n+1) \\
&\times
P_{n+j-l}^{(a,b+c+2l+1)}(1-2x_1)\
(1-x_1)^l P_l^{(b,c)}\left( {{1-x_1-2x_2}\over{1-x_1}}\right) \\
=&\, (x_1+x_2)^n P_n^{(a,b)}\left( {{x_2-x_1}\over{x_1+x_2}}\right)
\
P_j^{(a+b+2n+1,c)}(1-2(x_1+x_2)).
\endalign
$$
{\rm (ii)} The Hahn polynomials satisfy the following convolution
identity
$$
\align
\sum_{l=0}^{n+j}&
\left(j+n\atop l\right)
{(a+1)_{n+j-l}(b+1)_l (b+1)_n(c+1)_j(a+b+c+j+n+2)_l
\over (a+1)_n(c+1)_l(b+c+l+1)_l (b+c+2l+2)_{j+n-l}
(a+b+2n+2)_j }\\
&\times R_l(\la(n);b,c,-j-n-1,a+b+j+n+1) \\
&\times (l-s)_{n+j-l} \, Q_{n+j-l}(x_1;a,b+c+2l+1,s-l) \,
(x_1-s)_l\, Q_l(x_2;b,c,s-x_1) \\
=&\, (-x_1-x_2)_n\, Q_n(x_1;a,b,x_1+x_2) \,
(n-s)_j\, Q_j(x_1+x_2-n;a+b+2n+1,c,s-n)
\endalign
$$
with the notation \thetag{\vgldefHahnpols},
\thetag{\vgldefRacahpols}.
\endproclaim

\demo{Proof} The first result follows from the limit transition
of the continuous Hahn polynomials to the Jacobi polynomials.
Replace $x_i$ by $sx_i$ and let $s\to\infty$.
The second result follows by a similar substitution as in the
proof of Corollary \thmref{\corthmgenconvoidforMP}(ii).
\qed\enddemo

\demo{Remark \theoremname{\remcorthmgeneralconvoforcontinuousHahn}}
(i) Similar as in Remark \thmref{\remcorthmgenconvoidforMP}(iii)
we have that Corollary
\thmref{\corthmgeneralconvoforcontinuousHahn}(ii)
and Theorem \thmref{\thmgeneralconvoforcontinuousHahn}
can be obtained from each other by formal substitution.
Theorem \thmref{\thmgeneralconvoforcontinuousHahn}
can be obtained from Corollary
\thmref{\corthmgeneralconvoforcontinuousHahn}(i)
by a double application of the Mellin transform.
Moreover, Corollary
\thmref{\corthmgeneralconvoforcontinuousHahn}
can be proved as
Theorem \thmref{\thmgeneralconvoforcontinuousHahn}
by analysing the action of $X$ and $X_c$, cf.
Remark \thmref{\remcorthmgenconvoidforMP}(ii),
in the three-fold tensor product.

(ii) Dunkl \cite{\DunkPJM, Thm.~4.2, Prop.~5.4},
\cite{\DunkCJM, Thm.~1.7} has obtained Corollary
\thmref{\corthmgeneralconvoforcontinuousHahn},
and hence Theorem
\thmref{\thmgeneralconvoforcontinuousHahn},
by a different method. Dunkl \cite{\DunkPJM} obtains
the two-variable Hahn polynomials by
judiciously guessing solutions for a certain
difference equation arising from the representation theory
of the symmetric group. By symmetry considerations
there are more solutions of this type, and the
connection coefficients can be calculated in terms
of balanced ${}_4F_3$-series. The derivation in
this paper gives an intrinsic explanation for
the occurrence of the Racah polynomials as connection
coefficients. See also Dunkl \cite{\DunkPJM},
\cite{\DunkCJM} for the orthogonality
relations for these two-variable Hahn and Jacobi
polynomials for suitable restrictions on the parameters.
\enddemo

We do not obtain extensions of Corollary
\thmref{\corcorthmgenconvoidforMP} in this way.
For $k_1,k_2,k_3\to\infty$ in Theorem
\thmref{\thmgeneralconvoforcontinuousHahn}
we obtain the same result. This is also explained
by the fact that the Racah coefficients for
the Lie algebra ${\frak b}(1)$ are of the same
form as the Clebsch-Gordan coefficients,
cf. \cite{\VdJeug}.

\head\newsection The case $\Unc$\endhead

Let $\U$ be the complex unital
associative algebra generated by $A$, $B$, $C$, $D$ subject to the
relations
$$
AD=1=DA, \quad AB=qBA,\quad AC=q^{-1}CA,\quad
BC-CB = {{A^2-D^2}\over{q-q^{-1}}}.
\tag\eqname{\vgldefcommrelABCD}
$$
It is a Hopf algebra. We are only concerned with the
comultiplication, which is defined by
$$
\gathered
\De(A)=A\otimes A,\quad \De(B)=A\otimes B+B\otimes D,\\
\De(C) = A\otimes C+C\otimes D, \quad \De(D)=D\otimes D
\endgathered
\tag\eqname{\vgldefcomultiplication}
$$
on the level of generators and extended as an algebra
homomorphism. There are several possible $\ast$-structures
on $\U$, and we take
$$
A^\ast=A,\quad B^\ast=-C,\quad C^\ast=-B, \quad D^\ast=D,
$$
and the corresponding Hopf $\ast$-algebra is denoted by
$\Unc$.

The positive discrete series representations $\pi_k$ of $\Unc$
are unitary representations labelled by
$k>0$. They act in $\Hi$
and the action of the generators is given by
$$
\aligned
\pi_k(A)\, e^k_n &=q^{k+n}\, e^k_n, \\
\pi_k(C)\, e^k_n &= q^{1/2-k-n}
{{\sqrt{(1-q^{2n})(1-q^{4k+2n-2})}}\over{q-q^{-1}}}\, e^k_{n-1},\\
\pi_k(B) \, e^k_n &=
q^{-1/2-k-n}
{{\sqrt{(1-q^{2n+2})(1-q^{4k+2n})}}\over{q^{-1}-q}}\, e^k_{n+1}.
\endaligned
\tag\eqname{\vgldefqposdiscserreps}
$$
Note that $\pi_k(D)$ is an unbounded operator, but that
$\pi_k(A),\pi_k(B),\pi_k(C)\in {\Cal B}(\Hi)$. The operators that
we consider are bounded.

Recall that the tensor product of two representations
are defined by use of the comultiplication.
The tensor product of two positive discrete series representation
decomposes as for the Lie algebra ${\frak{su}}(1,1)$;
$$
\pi_{k_1}\otimes \pi_{k_2} \cong
\bigoplus_{j=0}^\infty \pi_{k_1+k_2+j}
\tag\eqname{\vglqtensorproddecomp}
$$
So there exists a unitary matrix mapping
the orthogonal basis $e_{n_1}^{k_1}\otimes e_{n_2}^{k_2}$ onto
$e_n^{k_1+k_2+j}$ intertwining the action of $\Unc$.
The matrix elements of this unitary mapping are
the Clebsch-Gordan coefficients;
$$
e_n^k = \sum_{n_1,n_2=0}^\infty C^{k_1,k_2,k}_{n_1,n_2, n}\,
e_{n_1}^{k_1}\otimes e_{n_2}^{k_2},
\tag\eqname{\vgldefqCGC}
$$
where $k=k_1+k_2+j$ for $j\in\Zp$. The sum is finite;
$n_1+n_2=n+j$. The Clebsch-Gordan coefficients are normalised
by $\langle e_0^k,e_0^{k_1}\otimes e_j^{k_2}\rangle >0$.

These results can be found in Burban and Klimyk \cite{\BurbK}
and Kalnins, Manocha and Miller \cite{\KalnMM}. See
Chari and Pressley \cite{\CharP} for general information on
quantised universal enveloping algebras.

\subhead\newsubsection
Clebsch-Gordan coefficients and orthogonal polynomials
\endsubhead
For this section we need the Askey-Wilson polynomials and the
Al-Salam and Chihara polynomials, which are a subclass of the
Askey-Wilson polynomials. The Askey-Wilson polynomial is
defined by
$$
p_m(\cos\th;a,b,c,d|q) = a^{-m} (ab,ac,ad;q)_m\, {}_4\vp_3
\left( {{q^{-m},abcdq^{m-1},ae^{i\th},ae^{-i\th}}\atop
{ab,\ ac,\ ad}}; q,q\right)
\tag\eqname{\vgldefAskeyWilsonpols}
$$
and it is symmetric in its parameters $a$, $b$, $c$ and $d$,
see \cite{\AskeW}.
The Al-Salam and Chihara polynomials are obtained by taking
$c=d=0$ in the Askey-Wilson polynomials;
$$
s_m(\cos\th;a,b|q) = p_m(\cos\th;a,b,0,0|q) = a^{-m} (ab;q)_m\,
{}_3\vp_2 \left( {{q^{-m},ae^{i\th},ae^{-i\th}}\atop{ab,\ 0}};
q,q\right).
\tag\eqname{\vgldefAlSalamChiharapols}
$$

By $dm(\cdot;a,b,c,d|q)$ we denote the normalised orthogonality
measure for the Askey-Wilson polynomials, which is absolutely
continuous on $[-1,1]$ and has at most a finite number of discrete
mass points outside $[-1,1]$. We put
$dm(\cdot;a,b|q)=dm(\cdot;a,b,0,0|q)$ for the normalised
orthogonality measure for the Al-Salam and Chihara polynomials.
Explicitly, let
$$
w(z) = {{(z^2,z^{-2};q)_\infty}\over
{(az,a/z,bz,b/z,cz,c/z,dz,d/z;q)_\infty}}
$$
and we use $w(z)=w(z;a,b,c,d|q)$ to stress the dependence on the
parameters when needed.
Let $a$, $b$, $c$ and $d$ be real, or, if complex,
appearing in conjugate pairs, and let all the pairwise products of
$a$, $b$, $c$ and $d$ not be greater or equal than $1$.
Then the Askey-Wilson polynomials
$p_n(x)= p_n(x;a,b,c,d| q)$ satisfy the orthogonality relations
$$
\aligned
{1\over {2\pi h_0}} \int_0^\pi & p_n(\cos\theta)p_m(\cos\theta)
w(e^{i\theta})\,d\theta
+{1\over{h_0}}\sum_k p_n(x_k)p_m(x_k)w_k = \delta_{n,m} h_n,\\
h_n &= {{ (1-q^{n-1}abcd)}\over{(1-q^{2n-1}abcd)}}
{{(q,ab,ac,ad,bc,bd,cd;q)_n}\over{(abcd;q)_n}} , \\
h_0 &= {{(abcd;q)_\infty}\over{(q,ab,ac,ad,bc,bd,cd;q)_\infty}}.
\endaligned
\tag\eqname{\vglorthoAskeyWilsonpols}
$$
The points $x_k$ are of the form ${1\over 2}(eq^k+e^{-1}q^{-k})$
for $e$ any of the parameters $a$, $b$, $c$ or $d$ with absolute
value greater than $1$; the sum is over $k\in\Zp$ such that
$\vert eq^k\vert >1$ and $w_k$ is the residue of $z\mapsto w(z)$ at
$z=eq^k$ minus the residue at $z=e^{-1}q^{-k}$. So
the normalised orthogonality measure
$dm(\cdot;a,b,c,d|q)$ can be read off from
\thetag{\vglorthoAskeyWilsonpols}, see Askey and Wilson
\cite{\AskeW} or \cite{\GaspR}.

Let $S_m(x;a,b|q)=s_m(x;a,b|q)/\sqrt{(q,ab;q)_m}$ denote the
orthonormal Al-Salam and Chihara polynomials, which satisfy
the three-term recurrence relation
$$
\aligned
2x\, S_n(x) &= a_{n+1}\, S_{n+1}(x) + q^n(a+b)\, S_n(x) +
a_n\, S_{n-1}(x),\\ a_n &= \sqrt{ (1-abq^{n-1})(1-q^n)}.
\endaligned
\tag\eqname{\vgldrietermAlSC}
$$

We now define
$$
Y_s = q^{1/2}B-q^{-1/2}C + {{s^{-1}+s}\over{q^{-1}-q}}(A-D)
\in \Unc.
\tag\eqname{\vgldefYs}
$$
Then $Y_sA$ is a self-adjoint element in $\Unc$ for
$s\in\R\backslash\{ 0\}$, or
$s\in\T$. $Y_s$ is twisted primitive element, i.e.
$\De(Y_s)=A\otimes Y_s+Y_s\otimes D$,
meaning that $Y_s$ is very much
like a Lie algebra element.

We also use the
notation $\mu(x)=(x+x^{-1})/2=\mu(x^{-1})$ for $x\not= 0$
in this section.

\proclaim{Proposition \theoremname{\propspectralYsA}}
Let $\Lambda\colon \Hi \to L^2(\R,dm(\cdot;q^{2k}s,q^{2k}/s|q^2))$
be the unitary mapping defined by
$\Lambda \colon e^k_n \mapsto S_n (\cdot;q^{2k}s,q^{2k}/s|q^2)$,
then $\Lambda$ intertwines $\pi_k(Y_sA)$ acting in
$\Hi$
with $2(M_x-\mu(s))/(q^{-1}-q)$.
\endproclaim

\demo{Proof} The bounded self-adjoint operator
$\pi_k(Y_sA)$ is a Jacobi matrix
by \thetag{\vgldefYs} and \thetag{\vgldefqposdiscserreps},
and the result follows upon comparing
with the three-term recurrence \thetag{\vgldrietermAlSC}
for the Al-Salam and Chihara polynomials as in
\S 2. 
\qed\enddemo

Proposition \thmref{\propspectralYsA} says that
$v^k(x) = \sum_{n=0}^\infty S_n( \mu(x); q^{2k}s,
q^{2k}/s|q^2)\, e_n^k$
is a generalised eigenvector of the self-adjoint operator
$\pi_k(Y_sA)$ for the eigenvalue
$$
\la_x = {{ x+x^{-1} - s - s^{-1}}\over{q^{-1}-q}} =
2{{\mu(x)-\mu(s)}\over{q^{-1}-q}}.
$$

Due to the fact that the comultiplication on $\Unc$
is less simple than for the Lie algebra ${\frak{su}}(1,1)$ it takes
a little more effort to determine the action of $Y_sA$ in
$\pi_{k_1}\otimes\pi_{k_2}$. The result can still be phrased
using orthogonal polynomials in two variables.

\proclaim{Proposition \theoremname{\propspectralYsAintensor}}
Define $\Upsilon\colon \Hi\otimes\Hi
\to L^2(\R^2,dm(x_1,x_2))$,
where
$$
dm(x_1,x_2) = dm(x_1;q^{2k_1}w_2,q^{2k_1}/w_2|q^2)\,
dm(x_2;q^{2k_2}s,q^{2k_2}/s|q^2), \qquad x_2=\mu(w_2),
$$
by
$$
\Upsilon \colon
e^{k_1}_{n_1} \otimes e^{k_2}_{n_2}
\mapsto S_{n_1} (x_1;q^{2k_1}w_2,q^{2k_1}/w_2|q^2)\,
S_{n_2}(x_2;q^{2k_2}s,q^{2k_2}/s|q^2)
$$
then $\Upsilon$ is a unitary mapping intertwining
$\pi_{k_1}\otimes\pi_{k_2}(\De(Y_sA))$
with $2(M_{x_1}-\mu(s))/(q^{-1}-q)$ in $L^2(dm(x_1,x_2))$.
\endproclaim

Note that $\Upsilon(e^{k_1}_{n_1} \otimes e^{k_2}_{n_2})$
forms a set of orthogonal
polynomials in two variables $x_1$ and $x_2$ for
$L^2(\R^2,dm(x_1,x_2))$, since the Al-Salam and Chihara polynomial
is symmetric in its parameters.

Proposition \thmref{\propspectralYsAintensor} states that
the vector
$$
\multline
w(x_1;x_2) = \sum_{n_1=0}^\infty S_{n_1}
(\mu(x_1);q^{2k_1}x_2,q^{2k_1}/x_2|q^2)\,
e^{k_1}_{n_1} \otimes v^{k_2}(x_2) \\
= \sum_{n_1,n_2=0}^\infty
S_{n_1} (\mu(x_1);q^{2k_1}x_2,q^{2k_1}/x_2|q^2)\,
S_{n_2}(\mu(x_2);q^{2k_2}s,q^{2k_2}/s|q^2)\,
e^{k_1}_{n_1} \otimes e^{k_2}_{n_2}
\endmultline
$$
is a generalised eigenvector of
$\pi_{k_1}\otimes\pi_{k_2}(\De(Y_sA))$
for the eigenvalue $\la_{x_1}$. This last observation is
essentially the way to obtain Proposition
\thmref{\propspectralYsAintensor}, since
$\De(Y_sA)=A^2\otimes Y_sA+Y_sA\otimes 1$ acts as a three-term
recurrence operator in $e^{k_1}_{n_1} \otimes v^{k_2}(x_2)$.

\demo{Proof} We use $\De(Y_sA)=A^2\otimes Y_sA+Y_sA\otimes 1$
and Proposition \thmref{\propspectralYsA} to
define for fixed $x_2$ the map
$\Lambda_0\colon\Hi\otimes\Hi\to \Hi$ by
$$
\Lambda_0 \colon e^{k_1}_{n_1} \otimes e^{k_2}_{n_2}
\mapsto S_{n_2}(x_2;q^{2k_2}s,q^{2k_2}/s|q^2)\, e^{k_1}_{n_1}
$$
to obtain the recurrence in $n_1$
$$
\multline
\Lambda_0 \Bigl( (q^{-1}-q)\bigl(\pi_{k_1}\otimes\pi_{k_2}
(\Delta(Y_sA)\bigr) + s +s^{-1}\Bigr)\,e^{k_1}_{n_1} \otimes
e^{k_2}_{n_2}  = \\ S_{n_2}(x_2;q^{2k_2}s,q^{2k_2}/s|q^2)
\Bigl( q^{2n_1} \bigl( (s+s^{-1})q^{2k_1} + \la_{w_2} q^{2k_1}
(q^{-1}-q)\bigr)\, e^{k_1}_{n_1} \\
 +\sqrt{ (1-q^{2n_1+2})(1-q^{4k_1+2n})}\, e_{n_1+1}^{k_1}
+ \sqrt{ (1-q^{2n_1})(1-q^{4k_1+2n_1-1})}\, e_{n_1-1}^{k_1}\Bigr)
\endmultline
$$
Use the explicit expression for $\la_{w_2}$ and the three-term
recurrence relation \thetag{\vgldrietermAlSC} to obtain the result.
\qed\enddemo

We now calculate the action of $\Upsilon\, e^k_n$, which yields
another set of orthonormal polynomials for $L^2(\R^2,dm(x_1,x_2))$.

\proclaim{Proposition \theoremname{\propactioUpsionlonekn}}
Let $k=k_1+k_2+j$ for $j\in\Zp$, and $x_1=\mu(w_1)$ then
$$
\align
\bigl( \Upsilon\, e^k_n\bigr)(x_1,x_2) &=
S_n(x_1;q^{2k}s,q^{2k}/s|q^2)
\, \bigl( \Upsilon\, e^k_0\bigr)(x_1,x_2), \\
\bigl( \Upsilon\, e^k_0\bigr)(x_1,x_2) &= C\,
p_j(x_2;q^{2k_1}w_1, q^{2k_1}/w_1,q^{2k_2}s,q^{2k_2}/s|q^2),\\
C^{-1} = \bigl((C_j(k_1,k_2)\bigr)^{-1}
&= \sqrt{ (q^2,q^{4k_1},q^{4k_2},q^{4k_1+4k_2+2j-2};q^2)_j}.
\endalign
$$
\endproclaim

\demo{Proof} By Proposition \thmref{\propspectralYsAintensor},
\thetag{\vglqtensorproddecomp} and
$$
2{{x_1-\mu(s)}\over{q^{-1}-q}} \, \Upsilon e_n^k(x_1,x_2) =
\Upsilon (\pi_k(Y_sA)\, e_n^k) (x_1,x_2)
$$
we obtain the three-term recurrence relation as
in Proposition \thmref{\propspectralYsA}, but with different
initial conditions. Hence, the first statement follows.

Since $\Upsilon$ is unitary we have the orthogonality relations
$\de_{nm}\de_{kl} = \langle \Upsilon e^k_n,\Upsilon e^l_m\rangle
=$
$$
\int S_n(x_1;q^{2k}s,q^{2k}/s|q^2)S_m(x_1;q^{2l}s,q^{2l}/s|q^2)
\int \Upsilon e^k_0(x_1,x_2) \Upsilon e^l_0(x_1,x_2)\, dm(x_1,x_2),
$$
by our first observation. As in the proof of Proposition
\thmref{\propspectdiffbasisDeltaXphi} we conclude
$\Upsilon e^k_0(x_1,x_2)=p_j(x_2)$ is polynomial of degree $j$,
$k=k_1+k_2+j$, in $x_2$ satisfying the orthogonality relations
$$
\int_{x_2} p_j(x_2)p_i(x_2)\, dm(x_1,x_2) = \de_{ij}
dm(x_1;q^{2k}s,q^{2k}/s|q^2)
$$
as measures with respect to functions in the variable $x_1$.

We now assume for ease of presentation
that $dm(x_1,x_2)$ is absolutely continuous. The general case can be
proved similarly, or it can be obtained by analytic continuation
with respect to $s$.
The measure is absolutely continuous for
$q^{2k_2}<|s|<q^{-2k_2}$, since $k_1,k_2>0$.
Put $x_1=\cos\th$, $x_2=\cos\psi$, then we obtain the explicit
expression \thetag{\vglorthoAskeyWilsonpols}
for the orthogonality measure;
$$
\multline
{1\over{2\pi}}\int_0^\pi p_i(\cos\psi)p_j(\cos\psi)
{{ (e^{\pm 2i\psi}, e^{\pm 2i\th};q^2)_\infty}\over
{(q^{2k_2}se^{\pm i\psi}, q^{2k_2}e^{\pm i\psi}/s,
q^{2k_1}e^{\pm i\psi \pm i\th};q^2)_\infty}} d\psi
= \\ \de_{ij}
{{ (q^{4k_1+4k_2+4j};q^2)_\infty}\over{
(q^2,q^{4k_1},q^{4k_2};q^2)_\infty}}
{{ (e^{\pm 2i\th};q^2)_\infty}\over{
(q^{2k_1+2k_2+2j}se^{\pm i\th}, q^{2k_1+2k_2+2j}
e^{\pm i\th}/s;q^2)_\infty}}
\endmultline
$$
for almost all $\th$. The $\pm$-signs means that we take all
possible combinations in the infinite $q$-shifted factorials.
Cancelling the $(e^{\pm 2i\th};q^2)_\infty$ on both sides and
comparing the result with \thetag{\vglorthoAskeyWilsonpols} we
see that $p_j$ is a multiple of
$p_j(\cdot;q^{2k_1}e^{i\th}, q^{2k_1}e^{-i\th},
q^{2k_2}s,q^{2k_2}/s|q^2)$.
The constant in front follows up to a sign by comparing the
squared norms. As in the proof of Proposition
\thmref{\propspectdiffbasisDeltaXphi} the sign of $C$ follows
from the normalisation of the Clebsch-Gordan coefficients, and now
we obtain $C>0$.
\qed\enddemo

So we obtain a second set of orthonormal polynomials for
$L^2(\R^2,dm(x_1,x_2))$ in terms of Al-Salam and Chihara
polynomials and Askey-Wilson polynomials.

The convolution formula for the Al-Salam and Chihara polynomials
is obtained by applying $\Upsilon$ to
\thetag{\vgldefqCGC} using the results of Propositions
\thmref{\propspectralYsAintensor} and
\thmref{\propactioUpsionlonekn}.
The results holds as an identity in a weighted $L^2$-space, but
since it is a polynomial identity it holds for all $x_1$, $x_2$;
with $x_1=\mu(w_1)$, $x_2=\mu(w_2)$, and $k=k_1+k_2+j$
$$
\multline
\sum_{n_1+n_2=n+j} C^{k_1,k_2,k}_{n_1,n_2,n}
\, S_{n_1} (x_1;q^{2k_1}w_2,q^{2k_1}/w_2|q^2)\,
S_{n_2}(x_2;q^{2k_2}s,q^{2k_2}/s|q^2) = \\
{{S_n(x_1;q^{2k}s,q^{2k}/s|q^2)
p_j(x_2;q^{2k_1}w_1, q^{2k_1}/w_1,q^{2k_2}s,q^{2k_2}/s|q^2)}\over{
 \sqrt{ (q^2,q^{4k_1},q^{4k_2},q^{4k_1+4k_2+2j-2};q^2)_j}}}.
\endmultline
\tag\eqname{\vgllemconvthmone}
$$

We have not yet calculated the Clebsch-Gordan
coefficients explicitly, but we can now use
\thetag{\vgllemconvthmone} to
determine $C^{k_1,k_2,k}_{n_1,n_2,n}$ by specialising to
a generating function for the Clebsch-Gordan coefficients.
The result is phrased in terms of $q$-Hahn polynomials, which are
defined as follows;
$$
Q_n(q^{-x};a,b,N;q) =
{}_3\vp_2 \left( {{q^{-n}, q^{-x},abq^{n+1}}\atop{
aq,\ q^{-N}}};q,q\right).
$$
See e.g. \cite{\KalnMM}
for other derivations of the following lemma.

\proclaim{Lemma \theoremname{\lemexplicitCGCgenfunct}}
With $n_1+n_2=n+j$ we get
$$
C^{k_1,k_2,k_1+k_2+j}_{n_1,n_2,n} =
C\ Q_j(q^{-2n_1};q^{4k_1-2},q^{4k_2-2},n+j;q^2),
$$
with the constant $C$ given by
$$
{{ q^{2k_1(n-n_1)}
(q^2;q^2)_{n+j}
\sqrt{(q^{4k_1};q^2)_{n_1} (q^{4k_2};q^2)_{n_2} (q^{4k_1};q^2)_j}}
\over {\sqrt{ (q^2;q^2)_n(q^2;q^2)_{n_1}(q^2;q^2)_{n_2}(q^2;q^2)_j
(q^{4k_1+4k_2+4j};q^2)_n (q^{4k_2};q^2)_j
(q^{4k_1+4k_2+2j-2};q^2)_j}}} .
$$
\endproclaim

\demo{Proof} Observe that $C^{k_1,k_2,k}_{n_1,n_2,n}$ is
independent of $s$, $x_1=\mu(w_1)$ and $x_2=\mu(w_2)$.
Specialise $w_2=q^{2k_2}s$ and $w_1=q^{2k_1}/w_2=q^{2k_1-2k_2}/s$,
then the Al-Salam and Chihara polynomials in the summand on the
left hand side of \thetag{\vgllemconvthmone}
can be evaluated explicitly, since the ${}_3\vp_2$-series
reduces to $1$. For this choice the Askey-Wilson polynomial
on the right hand side can also be evaluated explicitly,
and we obtain the generating function for the
Clebsch-Gordan coefficients
$$
\multline
\sum_{n_1+n_2=n+j} C^{k_1,k_2,k}_{n_1,n_2,n}
\, q^{2n_1(k_2-k_1)-2n_2k_2} s^{n_1-n_2}
{{ \sqrt{ (q^{4k_1};q^2)_{n_1}(q^{4k_2};q^2)_{n_2}}}\over
 { \sqrt{ (q^2;q^2)_{n_1}(q^2;q^2)_{n_2}}}} = \\
{{q^{-2jk_2 -2n(k_1+k_2+j)} s^{n-j}
(q^{4k_1},q^{4k_2}, q^{4k_2}s^2;q^2)_j }\over{
 \sqrt{ (q^2,q^{4k_1},q^{4k_2},q^{4k_1+4k_2+2j-2};q^2)_j}}}
\sqrt{ {{ (q^{4k_1+4k_2+4j};q^2)_n}\over{(q^2;q^2)_n}}}\\ \times
\ {}_3\vp_2 \left( {{q^{-2n},q^{4k_2+2j}, q^{4k_1+2j}/s^2}\atop{
q^{4k_1+4k_2+4j},\ 0}};q^2,q^2\right).
\endmultline
$$
This determines $C^{k_1,k_2,k}_{n_1,n_2,n}$, but it takes some
work to find the expression in terms of $q$-Hahn polynomials.
First, take $n_1$ as the summation parameter in the sum and
multiply both sides by $s^{n+j}$ to find that both sides are
polynomials of degree $n+j$ in $s^2$. Apply \cite{\GaspR, (III.6)}
to rewrite the ${}_3\vp_2$-series as a polynomial in $s^2$
and the $q$-binomial
theorem \cite{\GaspR, (II.3)} to write $(q^{4k_2}s^2;q^2)_j$ as
a polynomial in $s^2$. Comparing next the coefficients
on both sides gives an expression for the Clebsch-Gordan
coefficients as a terminating ${}_3\vp_2$-series. To put it into the
required form in terms of $q$-Hahn polynomials, we need to apply
some transformations for ${}_3\vp_2$-series, namely
\cite{\GaspR, (III.13), (III.11)}. The constant follows by
a straightforward calculation.
\qed\enddemo

Combining Lemma \thmref{\lemexplicitCGCgenfunct} with the
unitarity of the intertwining operator consisting of the
Clebsch-Gordan coefficients results in the orthogonality
relations for the $q$-Hahn and dual $q$-Hahn
polynomials, cf. \cite{\GaspR, \S 7.2}.

We now have all ingredients to rewrite \thetag{\vgllemconvthmone}.
Simplifying proves the following theorem.

\proclaim{Theorem \theoremname{\thmgenconvoforAlSCpols}}
With the notation \thetag{\vgldefAlSalamChiharapols}
and \thetag{\vgldefAskeyWilsonpols} for
the Al-Salam and Chihara polynomials and Askey-Wilson
polynomials and
$x_1=\mu(w_1)$, $x_2=\mu(w_2)$, $n,j\in\Zp$, $k_1,k_2>0$ we have
$$
\multline
(q^{4k_1};q^2)_j \sum_{l=0}^{n+j}
q^{2k_1(n-l)}
\left[ {{n+j}\atop l}\right]_{q^2}\
Q_j(q^{-2l};q^{4k_1-2}, q^{4k_2-2},n+j;q^2)\\
\times
s_l(x_1;q^{2k_1}w_2,q^{2k_1}/w_2|q^2)\,
s_{n+j-l}(x_2;q^{2k_2}s,q^{2k_2}/s|q^2) =
\\ s_n(x_1;q^{2k_1+2k_2+2j}s,q^{2k_1+2k_2+2j}/s|q^2)\,
p_j(x_2;q^{2k_1}w_1,q^{2k_1}/w_1,q^{2k_2}s,q^{2k_2}/s|q^2).
\endmultline
$$
\endproclaim

\demo{Remark \theoremname{\remthmgenconvoforAlSCpols}}(i)
Theorem \thmref{\thmgenconvoforAlSCpols} is a connection
coefficient formula for orthogonal polynomials in two
variables, orthogonal for the same measure, where the
connection coefficients are given by the
$q$-Hahn polynomials. The dual connection coefficient problem
follows from the orthogonality for the Clebsch-Gordan coefficients
or, equivalently, from the orthogonality
relations for the dual $q$-Hahn polynomials.

(ii) The case $j=0$ gives a simple convolution
property for the Al-Salam and Chihara polynomials, since
the $q$-Hahn and the Askey-Wilson polynomial reduce to $1$.
This
was the motivation for Al-Salam and Chihara \cite{\AlSaC} to
introduce the Al-Salam and Chihara polynomials as the
most general set of orthogonal polynomials still
satisfying a convolution property, see also
Al-Salam \cite{\AlSa, \S 8}.
The case $n=0$ is also of interest, since then the
$q$-Hahn polynomial can be evaluated and the Al-Salam and
Chihara polynomial on the right hand side reduces to $1$.
In both cases we have a free parameter in the sum.

(iii) Formally, in the representation space $\Hi\otimes\Hi$ we have
two bases of (generalised) eigenvectors for the action of $Y_sA$,
namely $v^k(x)$ and $w(x_1;x_2)$.
They are connected by Clebsch-Gordan coefficients, which are
now expressible as Askey-Wilson polynomials;
$$
w(x_1;x_2) = \sum_{j=0}^\infty {{
p_j(\mu(x_2);q^{2k_1}x_1, q^{2k_1}/x_1,q^{2k_2}s,q^{2k_2}/s|q^2)}
\over{ \sqrt{ (q^2,q^{4k_1},q^{4k_2},q^{4k_1+4k_2+2j-2};q^2)_j}}}
\, v^{k_1+k_2+j}(x_1).
$$
The dual Clebsch-Gordan coefficient relation follows by
integrating against the appropriate orthogonality measure for the
Askey-Wilson polynomials.
\enddemo

\subhead\newsubsection Racah coefficients and orthogonal
polynomials \endsubhead
In the tensor product of three positive discrete
series representations $\pi_{k_1} \otimes
\pi_{k_2} \otimes \pi_{k_3}$ of $\Unc$
we have the same orthogonal bases as in \S 3.2 
and we use the same notation as in
\thetag{\vgldeforthobasisdrieeen}--
\thetag{\vgldeforthobasisdrievier}.
Similarly, we now have an intertwiner for the $\Unc$-action
in terms of $q$-Racah coefficients;
$$
e^{((k_1k_2)k_{12}k_3)k}_n = \sum_{k_{23}}
U^{k_1,k_2,k_{12}}_{k_3,k,k_{23}} \ e^{(k_1(k_2k_3)k_{23})k}_n .
\tag\eqname{\vgldefqracahcoeff}
$$
Again the constraints
\thetag{\vglconventiosumvariables} hold.

\proclaim{Proposition \theoremname{\propThetaforUqsu}}
Define $\Theta\colon \Hi\otimes \Hi \otimes \Hi
\rightarrow L^2(\R^3,dm(x_1,x_2,x_3))$ by
$$
\multline
\Theta\bigl(e^{k_1}_{n_1} \otimes e^{k_2}_{n_2} \otimes
e^{k_3}_{n_3} \bigr)(x_1,x_2,x_3)\\
= S_{n_1}(x_1;q^{2k_1}w_2,q^{2k_1}/w_2|q^2)\,
S_{n_2}(x_2;q^{2k_2}w_3,q^{2k_2}/w_3|q^2)\,
S_{n_3}(x_3;q^{2k_3}s,q^{2k_3}/s|q^2),
\endmultline
$$
with the measure
$dm(x_1,x_2,x_3)$ given by
$$
dm(x_1;q^{2k_1}w_2,q^{2k_1}/w_2|q^2)\,
dm(x_2;q^{2k_2}w_3,q^{2k_2}/w_3|q^2)\,
dm(x_3;q^{2k_3}s,q^{2k_3}/s|q^2),
$$
where $x_i=\mu(w_i)$.
Then $\Theta$ is a unitary map intertwining $\pi_{k_1} \otimes
\pi_{k_2} \otimes \pi_{k_3}((1 \otimes \Delta)(\Delta(Y_s A))$
with $2(M_{x_1}-\mu(s))/(q^{-1}-q)$.
\endproclaim

\demo{Proof}
Observe that $(1 \otimes \Delta)(\Delta(Y_s A)) =
A^2\otimes \Delta(Y_s A) + Y_s A \otimes \Delta(1)$.
The proof now proceeds as the proof of
Proposition \thmref{\propspectralYsAintensor}.
\qed\enddemo

\proclaim{Proposition \theoremname{\propqThetaonbases}}
{\rm (i)} The following equality holds with $x_i=\mu(w_i)$;
$$
\gather
\Theta \bigl(e^{((k_1k_2)k_{12}k_3)k}_n\bigr)(x_1,x_2,x_3) =
C_j(k_{12},k_3) C_{j_{12}}(k_1,k_2)\, S_n(x_1;q^{2k}s,q^{2k}/s|q^2)
\\ \times
p_j(x_3;q^{2k_{12}}w_1,q^{2k_{12}}/w_1,q^{2k_{3}}s,q^{2k_{3}}/s|q^2)
p_{j_{12}}(x_2;q^{2k_1}w_1,q^{2k_1}/w_1,
q^{2k_{2}}w_3,q^{2k_{2}}/w_3|q^2).
\endgather
$$
{\rm (ii)} The following equality holds with $x_i=\mu(w_i)$;
$$
\gather
\Theta \bigl(e^{(k_1(k_2k_3)k_{23})k}_n\bigr) (x_1,x_2,x_3) =
C_{j'}(k_1,k_{23}) C_{j_{23}}(k_2,k_3)\,
S_n(x_1;q^{2k}s,q^{2k}/s|q^2)\\ \times
p_{j'}(x_2;q^{2k_1}w_1,q^{2k_1}/w_1,q^{2k_{23}}s,q^{2k_{23}}/s|q^2)
p_{j_{23}}(x_3;q^{2k_2}w_2,q^{2k_2}/w_2,
q^{2k_{3}}s,q^{2k_{23}}/s|q^2).
\endgather
$$
The constant $C_j(k_1,k_2)$ is defined in Proposition
\thmref{\propactioUpsionlonekn}.
\endproclaim

The proof of Proposition \thmref{\propqThetaonbases}
is slightly more complicated than the proof of its
counterpart Proposition
\thmref{\propspectdiffbasisDeltaDeltaXphi} due to the
fact that we do not have a nice factorisation for
$\Theta$ as in Remark
\thmref{\rempropspceDeltaDeltaXphi}. This is a
consequence of the non-cocommutativity of the
comultiplication for $\Unc$.

Note that the occurrence of $S_n(x_1;q^{2k}s,q^{2k}/s|q^2)$
on the right hand side corresponds to the intertwining
property of the Racah coefficients as in the proof of the
first statement of Proposition \thmref{\propactioUpsionlonekn}.

\demo{Proof} The proof of (i) and (ii) is similar. To prove
(i) we use \thetag{\vgldeforthobasisdrietwee} (for the
$\Unc$-setting) and Proposition
\thmref{\propThetaforUqsu} to find
$$
\align
&\Theta \bigl(e^{((k_1k_2)k_{12}k_3)k}_n\bigr)(x_1,x_2,x_3) =
\sum_{n_{12}+n_3=n+j} C^{k_{12},k_3,k}_{n_{12},n_3,n}
\, S_{n_3}(x_3;q^{2k_3}s,q^{2k_3}/s|q^2) \\
&\quad \times \sum_{n_1+n_2=n_{12}+j_{12}}
C^{k_1,k_2,k_{12}}_{n_1,n_2,n_{12}}\,
S_{n_1}(x_1;q^{2k_1}w_2,q^{2k_1}/w_2|q^2)\,
S_{n_2}(x_2;q^{2k_2}w_3,q^{2k_2}/w_3|q^2) \\
=& C_{j_{12}}(k_1,k_2)
p_{j_{12}}(x_2;q^{2k_1}w_1,q^{2k_1}/w_1,
q^{2k_{2}}w_3,q^{2k_{2}}/w_3|q^2)\\ &\quad\times
\sum_{n_{12}+n_3=n+j} C^{k_{12},k_3,k}_{n_{12},n_3,n}
\, S_{n_3}(x_3;q^{2k_3}s,q^{2k_3}/s|q^2)
\, S_{n_{12}}(x_1;q^{2k_{12}}w_3,q^{2k_{12}}/w_3,|q^2)
\endalign
$$
by \thetag{\vgllemconvthmone}. The last sum can be evaluated
by another application of \thetag{\vgllemconvthmone}
leading to the result.
\qed\enddemo

The $n$-dependence in the right hand sides of Proposition
\thmref{\propqThetaonbases} is the same, so we obtain the
Wigner-Eckhart theorem for the $\Unc$-setting by applying
$\Theta$ to \thetag{\vgldefqracahcoeff}.
So we can restrict
to $n=0$ in \thetag{\vgldefqracahcoeff} before applying
$\Theta$ without loss of generality, and we obtain
the following polynomial identity in $x_2$ and $x_3$
with $w_1$ as a parameter;
$$
\multline
C_j(k_{12},k_3) C_{j_{12}}(k_1,k_2)\,
p_{j}(x_3;q^{2k_{12}}w_1,q^{2k_{12}}/w_1,
q^{2k_{3}}s,q^{2k_{3}}/s|q^2)\\ \times
p_{j_{12}}(x_2;q^{2k_1}w_1,q^{2k_1}/w_1,q^{2k_{2}}w_3,
q^{2k_{2}}/w_3|q^2)
= \sum_{j_{23}=0}^{j_{12}+j} U^{k_1,k_2,k_{12}}_{k_3,k,k_{23}}
\, C_{j'}(k_1,k_{23})C_{j_{23}}(k_2,k_3)\\
\times
p_{j'}(x_2;q^{2k_1}w_1,q^{2k_1}/w_1,q^{2k_{23}}s,
q^{2k_{23}}/s|q^2)
\, p_{j_{23}}(x_3;q^{2k_2}w_2,q^{2k_2}/w_2,
q^{2k_{3}}s,q^{2k_{23}}/s|q^2).
\endmultline
\tag\eqname{\vglqconvoRCone}
$$
Again we can use \thetag{\vglqconvoRCone} in two ways. Firstly,
we specialise to a suitable formula from which the
Racah coefficients can be determined explicitly.
Secondly, with the explicit expression for the
Racah coefficients we derive a
convolution identity for the Askey-Wilson polynomials.

\proclaim{Proposition \theoremname{\propexplexprqRacahc}}
The Racah coefficients of \thetag{\vgldefqracahcoeff} are given by
$$
U^{k_1,k_2,k_{12}}_{k_3,k,k_{23}} = C\;
{}_4\varphi_3\left(
{q^{4k_1+4k_2+2j_{12}-2},q^{4k_2+4k_3+2j_{23}-2},q^{-2j_{12}},
q^{-2j_{23}} \atop q^{4k_2},q^{4k_1+4k_2+4k_3+2j+2j_{12}-2},
q^{-2j-2j_{12}}} ;q^2;q^2\right),
$$
with the constant $C$ given by
$$
\align
&{{(q^2,q^{4k_1},q^{4k_{23}},q^{4k_1+4k_{23}+2(j+j_{12}-j_{23})-2};
q^2)_{j+j_{12}-j_{23}}^{1/2}
(q^2,q^{4k_2},q^{4k_{3}},q^{4k_2+4k_{3}+2j_{23}-2};
q^2)_{j_{23}}^{1/2}} \over
{(q^2,q^{4k_{12}},q^{4k_{3}},q^{4k_{12}+4k_{3}+2j-2};q^2)_j^{1/2}
(q^2,q^{4k_1},q^{4k_{2}},q^{4k_1+4k_{2}+2j_{12}-2};
q^2)_{j_{12}}^{1/2}}} \\
&\quad \times q^{2k_2(j-j_{23})}\left[j+j_{12}\atop j_{23}
\right]_{q^2} {(q^{4k_3};q^2)_j(q^{4k_2};q^2)_{j_{12}}
(q^{4k_1+4k_2+4k_3+2j+2j_{12}-2};q^2)_{j_{23}}
\over (q^{4k_3},q^{4k_2+4k_3+2j_{23}-2};q^2)_{j_{23}}
(q^{4k_2+4k_3+4j_{23}}; q^2)_{j+j_{12}-j_{23}} }.
\endalign
$$
\endproclaim

\demo{Proof} Take $w_1=s=1$ and
$x_2=\mu(q^{2k_1})$ in \thetag{\vglqconvoRCone} to find
$$
\multline
p_j(x_3;q^{2k_{12}},q^{2k_{12}},q^{2k_3},q^{2k_3}|q^2)
(q^{2k_1+2k_2}w_3, q^{2k_1+2k_2}/w_3;q^2)_{j_{12}} \\
=\sum_{j_{23}}C_1 U^{k_1,k_2,k_{12}}_{k_3,k,k_{23}}
p_{j_{23}}(x_3;q^{2k_1+2k_2},q^{2k_2-2k_1},q^{2k_3},q^{2k_3}|q^2)
\endmultline
$$
for $C_1$ an explicit constant depending upon
$k_1,k_2,k_3,j_{12},j_{23}$ and $j$, since two Askey-Wilson
polynomials can be evaluated for this choice. So the Racah
coefficients occur as the coefficients when developing the
polynomial of degree $j+j_{12}$ on the left hand side
into Askey-Wilson polynomials. Hence the Racah coefficients
can be obtained from
$$
\multline
C_2 U^{k_1,k_2,k_{12}}_{k_3,k,k_{23}}=\\
\int p_{j_{23}}(x_3;q^{2k_1+2k_2},q^{2k_2-2k_1},
q^{2k_3},q^{2k_3}|q^2)
p_j(x_3;q^{2k_{12}},q^{2k_{12}},q^{2k_3},q^{2k_3}|q^2) \\
\times(q^{2k_1+2k_2}w_3, q^{2k_1+2k_2}/w_3;q^2)_{j_{12}}
\ dm(x_3;q^{2k_1+2k_2},q^{2k_2-2k_1},q^{2k_3},q^{2k_3}|q^2),
\endmultline
$$
for some known constant $C_2$. Now observe that
$$
\multline
(q^{2k_1+2k_2}w_3, q^{2k_1+2k_2}/w_3;q^2)_{j_{12}} \,
dm(x_3;q^{2k_1+2k_2},q^{2k_2-2k_1},q^{2k_3},q^{2k_3}|q^2)
\\ = C_3\, dm(x_3;q^{2k_{12}},q^{2k_2-2k_1},q^{2k_3},q^{2k_3}|q^2)
\endmultline
$$
for some constant known $C_3$ by \thetag{\vglconventiosumvariables}
and \thetag{\vglorthoAskeyWilsonpols}.
Thus the Racah coefficients can be obtained
by integration;
$$
\multline
C_4\, U^{k_1,k_2,k_{12}}_{k_3,k,k_{23}}=\int
p_{j_{23}}(x_3;q^{2k_1+2k_2},q^{2k_2-2k_1},q^{2k_3},q^{2k_3}|q^2)\\
\times p_j(x_3;q^{2k_{12}},q^{2k_{12}},q^{2k_3},q^{2k_3}|q^2)
\ dm(w_3;q^{2k_{12}},q^{2k_2-2k_1},q^{2k_3},q^{2k_3}|q^2),
\endmultline
\tag\eqname{\vgleenRacahasconcoefAWpols}
$$
with $C_4$ explicitly known. Observe that three out of four
of the parameters of each of the Askey-Wilson
polynomials in \thetag{\vgleenRacahasconcoefAWpols}
coincide with the parameters of the Askey-Wilson measure in
\thetag{\vgleenRacahasconcoefAWpols}. Use the connection
coefficient formula for Askey-Wilson polynomials
with one different parameter, cf. Askey and Wilson
\cite{\AskeW, (6.4-5)} or see \cite{\GaspR, (7.6.8-9)
with the right hand side of (7.6.9) multiplied
by $(q;q)_n$}, twice
to rewrite the Askey-Wilson polynomials in terms
of Askey-Wilson polynomials with the same parameters as
the Askey-Wilson measure in \thetag{\vgleenRacahasconcoefAWpols}.
By orthogonality the integration is then easily performed and we
are left with a single sum, which can be written as a
very-well-poised ${}_8\vp_7$-series. This can be transformed to a
balanced ${}_4\vp_3$-series by Watson's transformation
\cite{\GaspR, (III.17)}, and another application of
Sears's transformation \cite{\GaspR, (III.15)} gives the
form as in the statement of the proposition.
The constant follows from bookkeeping.
\qed\enddemo

Recall the $q$-Racah polynomials,
see \cite{\GaspR, \S 7.2}, \cite{\KoekS},
$$
R_n(\nu(x);\al,\be,\ga,\de;q) =
{}_4\vp_3\left( {{q^{-n}, \al\be q^{n+1}, q^{-x}, \ga\de q^{x+1}}
\atop {\al q,\ \be\de q,\ \ga q}};q,q\right)
\tag\eqname{\vgldefqRacahpols}
$$
with $\nu(x)= q^{-x}+\ga\de q^{x+1}$, one of the lower parameters
equals $q^{-N}$, $N\in\Zp$, and $0\leq n \leq N$.
The ${}_4\vp_3$-series in Proposition
\thmref{\propexplexprqRacahc} can be written in terms of
a $q$-Racah polynomial.

We can now rewrite \thetag{\vglqconvoRCone} to arrive at the
key result of this paper.
For convenience we
replace $q^2$ by $q$, $(a,b,c) = (q^{k_1},q^{k_2},q^{k_3})$, and we
relabel $w_1$, $j_{23}$ and $j_{12}$ by $t$, $l$ and $n$, and
finally replace $x_2$, $x_3$ by $x_1$, $x_2$. We obtain the
following $q$-analogue of Theorem
\thmref{\thmgeneralconvoforcontinuousHahn} and
Corollary \thmref{\corthmgeneralconvoforcontinuousHahn}.

\proclaim{Theorem \theoremname{\thmconvoforAskeyWilson}}
With $x_1=\mu(w_1)$, $x_2=\mu(w_2)$, $n,j\in\Zp$, we have
the convolution identity for the Askey-Wilson
polynomials
$$
\align
\sum_{l=0}^{n+j} & b^{j-l} \left[ j+n\atop l\right]_q
{(b^2;q)_n(a^2b^2c^2q^{j+n-1};q)_l(c^2;q)_j \over
(c^2,b^2c^2q^{l-1};q)_l (b^2c^2q^{2l};q)_{j+n-l} } \\ &\times
R_l(\nu(n);b^2/q,c^2/q,q^{-j-n-1},a^2b^2q^{n+j-1};q)\\ &\times
p_{j+n-l}(x_1;at,a/t,bcq^ls,bcq^l/s|q)\,
p_l(x_2;bw_1,b/w_1,cs,c/s|q)\\
=\, &p_n(x_1;at,a/t,bw_2,b/w_2|q)\,
p_j(x_2;abq^nt,abq^n/t,cs,c/s|q),
\endalign
$$
with the notation of \thetag{\vgldefAskeyWilsonpols},
\thetag{\vgldefqRacahpols}.
\endproclaim

\demo{Remark \theoremname{\remthmconvoforAskeyWilson}} (i)
Theorem \thmref{\thmgenconvoforAlSCpols} can be obtained
as a special case of Theorem \thmref{\thmconvoforAskeyWilson}
by letting $c\to 0$.

(ii) Theorem \thmref{\thmconvoforAskeyWilson} leads to
a kind of generating function for the $q$-Racah and $q$-Hahn
polynomials. Choosing $w_1=at$ and $w_2=cs$ in Theorem
\thmref{\thmconvoforAskeyWilson} reduces
all four Askey-Wilson polynomials to a single term. The remaining
free parameters $s$ and $t$ appear only in the combination $s/t$.
Replacing $bcs/(at)$ by $u$, and $(a^2,b^2,c^2)$ by $(\al,\be,\ga)$
gives:
$$
\multline
\sum_{l=0}^{n+j} \left[j+n\atop l\right]_q { (\al;q)_{j+n-l}
(\be;q)_n (\al\be\ga q^{j+n-1};q)_l \over (\al;q)_n
(\be\ga q^{l-1};q)_l (\be\ga q^{2l};q)_{j+n-l} } u^{j-l}
(\be\ga q^l/u;q)_{j+n-l} (u;q)_l \\
\times  R_l(\nu(n);\be/q,\ga/q,q^{-j-n-1},\al\be q^{n+j-1};q)
 = (\al q^n u;q)_j (\be/u;q)_n .
\endmultline
$$
For $\ga=0$ we obtain a similar identity for $q$-Hahn polynomials:
$$
\sum_{l=0}^{n+j} \left[j+n\atop l\right]_q { (\al;q)_{j+n-l}
(\be;q)_n \over (\al;q)_n }  Q_n(q^{-l};\be/q,\al/q,n+j;q)
u^{j-l} (u;q)_l = (\al q^n u;q)_j (\be/u;q)_n .
$$

(iii) Theorem \thmref{\thmconvoforAskeyWilson} gives the
connection coefficients for two sets of orthogonal
polynomials with respect to the absolutely continuous
measure
$$
{{(w_1^{\pm 2}, w_2^{\pm 2};q)_\infty}\over
{(taw_1^{\pm 1}, aw_1^{\pm 1}/t,
csw_2^{\pm 1}, cw_2^{\pm 1}/s;q)_\infty}}
{1\over{(bw_1^{\pm 1}w_2^{\pm 1};q)_\infty}} {{dw_1}\over{w_1}}
{{dw_2}\over{w_2}}
$$
on the torus $\T^2$ for $|t|^{-1}<|a|<|t|$, $|s|^{-1}<|c|<|s|$
and $|b|<1$. Here all possible signs for $\pm$ have to be used.
(Otherwise discrete masses at points and lines have to
be added.) This weight function is invariant under
simultaneously interchanging
$w_1$ with $w_2$, $t$ with $s$ and $a$ and $c$. This transforms
the orthogonal polynomials on the right hand side of
Theorem \thmref{\thmconvoforAskeyWilson} to the ones
occurring in the left hand side.

In particular, note that for $s=t$, $a=c$, the weight function
is invariant under the Weyl group for $B_2$, i.e. the group
generated by $(w_1,w_2)\mapsto (w_2,w_1)$ and
$(w_1,w_2)\mapsto (w_1,w_2^{-1})$. The corresponding Weyl group
invariant orthogonal polynomials are
$$
\multline
p_n(\mu(w_1);at,a/t,bw_2,b/w_2|q)\,
p_j(\mu(w_2);abq^nt,abq^n/t,at,a/t|q)
\\
+p_n(\mu(w_2);at,a/t,bw_1,b/w_1|q)\,
p_j(\mu(w_1);abq^nt,abq^n/t,at,a/t|q)
\endmultline
$$
for $n\geq j\geq 0$.
These orthogonal polynomials do not seem directly related to the
Koornwinder-Macdonald orthogonal polynomials associated with
root system $BC_2$, see \cite{\KoorContM}, although the
structure of the orthogonality measure is similar.
\enddemo

\head\newsection Linearisation coefficients for
Askey-Wilson polynomials\endhead

In the results of the previous section using $\Unc$,
especially Theorems \thmref{\thmgenconvoforAlSCpols}
and \thmref{\thmconvoforAskeyWilson}, we can use analytic
continuation with respect to the parameters involved
to find similar identities but with the Al-Salam
and Chihara polynomials and the Askey-Wilson polynomials
replaced by the dual $q$-Krawtchouk polynomials and
the $q$-Racah polynomials. These identities can be obtained
by the same procedure using $\Usu$ and its representation
theory instead of using $\Unc$.
In particular we can now give an interpretation
for $q$-Racah polynomials as Clebsch-Gordan coefficients
for $\Usu$. In this case we also have some
knowledge on the structure of the dual Hopf $\ast$-algebra
and this can be used to obtain a linearisation formula
for a two-parameter family of Askey-Wilson polynomials.
This is an application of the results of the previous
section.

We first recall $\Usu$ and its representation theory,
see e.g. \cite{\CharP}, \cite{\KoelAAM}, \cite{\KoorZSE}.
The Hopf algebra structure on $\Usu$ is the same as the
Hopf algebra structure on $\U$, cf.
\thetag{\vgldefcommrelABCD},
\thetag{\vgldefcomultiplication}. The $\ast$-operator
making $\Usu$ into a Hopf $\ast$-algebra is given by
$$
A^\ast=A,\quad B^\ast=C,\quad C^\ast=B, \quad D^\ast=D.
$$

There is precisely one irreducible unitary
$\Usu$-module $W_N$ of each dimension $N+1$
with highest weight vector $v_+$, i.e. $A\, v_+ = q^{N/2} v_+$,
$B\, v_+ = 0$. The corresponding representation is denoted
by $t^N$.
With respect to the standard orthonormal basis $e_n^N$,
$0\leq n\leq N$, the action of the generators is given by
$t^N(A)\, e^N_n = q^{n-N/2}\, e^N_d$ and
$$
\align
t^N(B)\, e^N_n &= {{q^{(1-N)/2}}\over{1-q^2}}
\sqrt{ (1-q^{2n+2})(1-q^{2N-2n})}\, e^N_{n+1}, \\
t^N(C)\, e^N_n &= {{q^{(1-N)/2}}\over{1-q^2}}
\sqrt{ (1-q^{2n})(1-q^{2N-2n+2}) }\, e^N_{n-1},
\endalign
$$
with the convention $e^N_{-1}=0=e^N_{N+1}$. So $e^N_N$ is the
highest weight vector. The representation $t^N$, considered as
a representation of $\U$ can be obtained from the
discrete series representation $\pi_k$ of
\thetag{\vgldefqposdiscserreps} by formally replacing
$k$ by $-N/2$.

The Clebsch-Gordan decomposition holds; as unitary $\Usu$-modules
$$
W_{N_1}\otimes W_{N_2} = \bigoplus_{j=0}^{\min(N_1,N_2)}
W_{N_1+N_2-2j}.
$$
The matrix coefficients of the intertwining operator give the
Clebsch-Gordan coefficients;
$$
e^N_n = \sum_{n_1,n_2} C^{N_1,N_2,N}_{n_1,n_2,n}\ e^{N_1}_{n_1}
\otimes e^{N_2}_{n_2}.
\tag\eqname{\vgldefCGCUqsutwee}
$$
Of course, the Clebsch-Gordan coefficient is
zero if $N\not= N_1+N_2-2j$ for $0\leq j\leq \min(N_1,N_2)$.
By considering the action of $A$ on both sides we see that the
Clebsch-Gordan coefficient is zero unless $n_1+n_2=n+j$, so that
the sum is actually a single sum. The Clebsch-Gordan coefficients
are normalised by
$\langle e^N_0,e^{N_1}_j\otimes e^{N_2}_0\rangle >0$
if $N=N_1+N_2-2j$, $0\leq j\leq \min(N_1,N_2)$.

We are particularly interested in the element
$$
X_p =q^{1/2}B + q^{-1/2}C -
{{p^{1/2}-p^{-1/2}}\over{q-q^{-1}}}(A-D) \in\U,\qquad p>0.
$$
Then $X_pA$ is self-adjoint and
$\De(X_pA)= A^2\otimes X_pA+X_pA\otimes 1$.
Koornwinder \cite{\KoorZSE} has shown that in each module $W_N$ the
action of $X_pA$ is completely diagonalisable. To formulate this
result we introduce the orthonormal dual $q$-Krawtchouk polynomials;
for $a>0$,
$$
r_n(q^{-x}-q^{x-N}/a;a,N;q) = (-1)^n
a^{n/2} q^{n(n-1)/4} \left[ {N\atop n}\right]_q^{1/2}
R_n(q^{-x}-q^{x-N}/a;a,N;q),
$$
for $N\in\Zp$ and $0\leq x,n\leq N$. The dual $q$-Krawtchouk
polynomials are defined by
$$
 R_n(q^{-x}-q^{x-N}/a;a,N;q) =
{}_3\vp_2 \left( {{q^{-n},q^{-x},-q^{x-N}/a}\atop
{q^{-N},0}};q,q\right).
$$
The corresponding three-term recurrence relation is
$$
\aligned
(q^{-x}-q^{x-N}/a)\, r_n &= A_n\, r_{n+1}
+ q^{n-N}(1-a^{-1})\, r_n + A_{n-1} \, r_{n-1}, \\
A_n & = a^{-1/2} q^{-N+ n/2} \sqrt{(1-q^{n+1})(1-q^{N-n})},
\endaligned
\tag\eqname{\vglthreetermdualqK}
$$
for $0\leq x,n\leq N$ and $r_n=r_n(q^{-x}-q^{x-N}/a;a,N;q)$.

\proclaim{Proposition \theoremname{\propspectrumXpAinWN}}
{\rm (\cite{\KoorZSE})} There exists
an orthogonal basis $\phi^N_f=\phi^N_f(p)$,
$0\leq f\leq N$,
of $W_N$ of eigenvectors of $t^N(X_pA)$ for the eigenvalue
$$
\la_f^N(p) = {{p^{1/2}q^{N-2f}-p^{-1/2}q^{2f-N}
+p^{-1/2}-p^{1/2}}\over{q^{-1}-q}}.
$$
Moreover, $\phi^N_f(p) = \sum_{n=0}^N
r_n(q^{-2f}-q^{2f-2N}/p;p,N;q^2)\, e_d^N$.
\endproclaim

Proposition \thmref{\propspectrumXpAinWN} is the analogue
of Proposition \thmref{\propspectralYsA}, and we could have
formulated it in a similar fashion using the finite discrete
orthogonality measure for the dual $q$-Krawtchouk polynomials.
Actually, replacing $e^{i\th}$, $k$ in
$S_n(\cos\th;q^{2k}s,q^{2k}/s|q^2)$
by $q^{-2f+N}/s$, $-N/2$ and next taking $s^2=-p^{-1}$ gives
$r_n(q^{-2f}-q^{2f-2N}/p;p,N;q^2)$.
The analogue of Proposition
\thmref{\propspectralYsAintensor} is the following.

\proclaim{Proposition \theoremname{\propeigvectXptensor}}
For $0\leq f_1\leq N_1$, $0\leq f_2\leq N_2$ define in
$W_{N_1}\otimes W_{N_2}$ the vector
$$
\phi^{N_1,N_2}_{f_1,f_2} = \sum_{n_1=0}^{N_1}
r_{n_1}(q^{-2f_1}-q^{2f_1-2N_1-2N_2+4f_2}/p;pq^{2N_2-4f_2},N_1;q^2)
\ e_{n_1}^{N_1}\otimes \phi_{f_2}^{N_2},
$$
then $t^{N_1}\otimes t^{N_2}\bigl( \De(X_pA)\bigr)
\, \phi^{N_1,N_2}_{f_1,f_2} =
\la_{f_1+f_2}^{N_1+N_2}(p)\, \phi^{N_1,N_2}_{f_1,f_2}$.
Moreover, $\phi^{N_1,N_2}_{f_1,f_2}$, $0\leq f_1\leq N_1$,
$0\leq f_2\leq N_2$, constitutes an orthogonal basis of
$W_{N_1}\otimes W_{N_2}$ of eigenvectors of $\Delta(X_pA)$.
\endproclaim

\demo{Proof} From $\De(X_pA)= A^2\otimes X_pA+X_pA\otimes 1$ it
follows that there is an eigenvector of the form
$\sum_{n_1=0}^{N_1} p_{n_1}\, e_{n_1}^{N_1}\otimes \phi_{f_2}^{N_2}$
by solving a three-term recurrence relation for the $p_{n_1}$.
Now \thetag{\vglthreetermdualqK} can be used to solve this.

There are $(N_1+1)(N_2+1)$ eigenvectors in $W_{N_1}\otimes W_{N_2}$
and $\langle \phi^{N_1,N_2}_{f_1,f_2},
\phi^{N_1,N_2}_{g_1,g_2}\rangle$ equals zero
if $f_2\not= g_2$ by Proposition \thmref{\propspectrumXpAinWN}
and it also equals zero if $f_1+f_2\not= g_1+g_2$ by the
self-adjointness of $X_pA$.
\qed\enddemo

The result of Proposition \thmref{\propeigvectXptensor} can be
obtained from Proposition \thmref{\propspectralYsAintensor} by
substituting $k_1$, $k_2$ by $-N_1/2$, $-N_2/2$ and $w_1$, $w_2$ by
$q^{N_1+N_2-2f_1-2f_2}/s$, $q^{N_2-2f_2}/s$ and $s^2$ by $-p^{-1}$.

\proclaim{Proposition \theoremname{\propeigvectDeXpeneNn}}
For $N=N_1+N_2-2j$, $0\leq j\leq\min(N_1,N_2)$ we have
$$
\langle \phi^{N_1,N_2}_{f_1,f_2}, e^N_n \rangle =
r_n(q^{2j-2f_1-2f_2}-q^{2f_1+2f_2-2j-2N}/p;p,N;q^2)
\langle \phi^{N_1,N_2}_{f_1,f_2}, e^N_0 \rangle,
$$
if $0\leq f_1+f_2-j\leq N$, and
$\langle \phi^{N_1,N_2}_{f_1,f_2}, e^N_n \rangle = 0$ otherwise.
If non-zero, then
$$
\multline
\langle \phi^{N_1,N_2}_{f_1,f_2}, e^N_0 \rangle =
\left[ {N_2\atop j}\right]_{q^2}^{1/2}
{{p^{j/2} q^{j(2N_1+N_2)} q^{-3j(j-1)/2}}\over
{\sqrt{ (q^{2N_1}, q^{2N_1+2N_2-2j+2};q^{-2})_j}}}
(-p^{-1}q^{2f_1+2f_2-2N_1-2N_2};q^2)_j
\\ \times (q^{-2f_1-2f_2};q^2)_j
\,  {}_4\vp_3
\left( {{q^{-2j}, q^{2j-2-2N_1-2N_2}, q^{-2f_2},
-p^{-1}q^{2f_2-2N_2}}\atop {q^{-2N_2}, q^{-2f_1-2f_2},
-p^{-1}q^{2f_1+2f_2-2N_1-2N_2}}};q^2,q^2\right).
\endmultline
$$
\endproclaim

Note that the ${}_4\vp_3$-series is balanced, and can be written
in terms of $q$-Racah polynomials \thetag{\vgldefqRacahpols}. The
${}_4\vp_3$-series in Proposition \thmref{\propeigvectDeXpeneNn}
equals
$$
R_j(q^{-2f_2}-p^{-1}q^{2f_2-2N_2};q^{-2N_2-2}, q^{-2N_1-2},
q^{-2f_1-2f_2-2}, -p^{-1}q^{2f_1+2f_2-2N_2};q^2).
$$
The proof of
Proposition \thmref{\propeigvectDeXpeneNn} is similar to the
proof of Proposition \thmref{\propactioUpsionlonekn}.
Proposition \thmref{\propeigvectDeXpeneNn} can also
be obtained from  Proposition \thmref{\propactioUpsionlonekn}
using the substitutions as indicated earlier.
It can also
be obtained by using $e^N_0=\sum C^{N_1,N_2,N}_{n_1,n_2,0}
\, e^{N_1}_{n_1}\otimes e^{N_2}_{n_2}$, Propositions
\thmref{\propeigvectXptensor} and \thmref{\propspectrumXpAinWN}
and the explicit value for the Clebsch-Gordan coefficients for
$n=0$, $N=N_1+N_2-2j$,
$$
C^{N_1,N_2,N}_{n_1,n_2,0} =
(-1)^{n_2} q^{n_2(N_2-j-1)}
\sqrt{ {{(q^{2n_1+2};q^2)_{n_2} (q^{2N_1-2n_1};q^{-2})_{n_2}}\over
{(q^2;q^2)_{n_2}(q^{2N_2};q^{-2})_{n_2}}}}
\sqrt{ {{(q^{2N_2};q^{-2})_j}\over{(q^{2N_1+N_2-2j+2};q^2)_j}}}.
$$
See e.g. \cite{\VileK, \S 14.3},
but this simple case can also be derived as follows.
Apply $C$ to both sides of \thetag{\vgldefCGCUqsutwee} to obtain
a three-term recurrence for the Clebsch-Gordan coefficients,
which reduces to a two-term recurrence for $n=0$. This can be
easily solved, with the initial condition following from
the unitarity and the normalisation, see \cite{\VdJeug} for
a similar derivation.
Then we have a sum involving the
product of two dual $q$-Krawtchouk polynomials. Upon inserting
the series representation we obtain a triple sum, and after
interchanging summations we can use the $q$-binomial theorem
and the $q$-Chu-Vandermonde sum, see \cite{\GaspR}, to obtain a
single ${}_4\vp_3$-series.

\demo{Remark \theoremname{\remconvosfromsutwo}}
With Proposition \thmref{\propeigvectDeXpeneNn} at hand
it is straightforward to calculate the $\Usu$-counterparts
of Theorems \thmref{\thmgenconvoforAlSCpols} and
\thmref{\thmconvoforAskeyWilson}. The Al-Salam and
Chihara, respectively Askey-Wilson, polynomials have to be
replaced by dual $q$-Krawtchouk, respectively $q$-Racah,
polynomials. The result can also be obtained from
Theorems \thmref{\thmgenconvoforAlSCpols} and
\thmref{\thmconvoforAskeyWilson} by substitution as
indicated earlier, so we do not give them explicitly.
These formulas give an alternative for the formulas
of Groza and Kachurik \cite{\GrozK}.
\enddemo

In the representation space $W_{N_1}\otimes W_{N_2}$ we have two
bases of eigenvectors for the action of $X_pA$, namely
$\phi_f^N$ and $\phi^{N_1,N_2}_{f_1,f_2}$, and the
corresponding Clebsch-Gordan coefficients are given
by Proposition \thmref{\propeigvectDeXpeneNn}, since,
with $N=N_1+N_2-2j$,
$$
\align
&\phi^{N_1,N_2}_{f_1,f_2} = \sum_{j=0}^{\min(N_1,N_2)}
\sum_{n=0}^N
\langle \phi^{N_1,N_2}_{f_1,f_2}, e^N_n\rangle \, e^N_n \\
&\ = \sum_{j=0}^{\min(N_1,N_2)}
\langle \phi^{N_1,N_2}_{f_1,f_2}, e^N_0\rangle
\sum_{n=0}^N r_n(q^{2j-2f_1-2f_2}-q^{2f_1+2f_2-2j-2N}/p;p,N;q^2)
\, e^N_n
\\ &\  = \sum_{j=0}^{\min(N_1,N_2)}
\langle \phi^{N_1,N_2}_{f_1,f_2}, e^N_0\rangle\, \phi^N_{f_1+f_2-j}.
\endalign
$$
Here we use the convention that $\phi^N_f=0$
for $f>N$ or $f<0$.
So introducing the notation
$$
\phi^{N_1,N_2}_{f_1,f_2} = \sum_{f,j} \,
C^{N_1,N_2,N}_{f_1,f_2,f}(p) \ \phi^N_f,
\tag\eqname{\vgldefgeneralisedCGCsutwo}
$$
we see that the Clebsch-Gordan coefficients are zero
unless $f_1+f_2=f+j$, and then
$C^{N_1,N_2,N}_{f_1,f_2,f}(p) =
\langle \phi^{N_1,N_2}_{f_1,f_2}, e^N_0\rangle$. So, by
Proposition \thmref{\propeigvectDeXpeneNn} we have proved
that the $q$-Racah polynomials occur as Clebsch-Gordan coefficients
for $\Usu$.

Using \thetag{\vgldefgeneralisedCGCsutwo} in a special
case we can obtain the linearisation coefficients for the
two parameter family of Askey-Wilson polynomials occurring
as spherical functions on the quantum $SU(2)$ group, cf.
\cite{\KoorZSE}. We consider odd-dimensional
representations; $N_1=2l_1$, $N_2=2l_2$, $l_1,l_2\in\Zp$.
Then the kernel of $t^{2l_1}(X_pA)$ is one dimensional and spanned
by $\phi^{2l_1}_{l_1}(p)$. Moreover, $\phi^{2l_1,2l_2}_{l_1,l_1}(p)=
\phi^{2l_1}_{l_1}(p)\otimes \phi^{2l_2}_{l_2}(p)$.
Next we consider matrix elements as linear functionals on
$\Usu$ to
find
$$
\aligned
&\sum_{(X)} \langle t^{2l_1}(X_{(1)})\phi^{2l_1}_{l_1}(p),
\phi^{2l_1}_{l_1}(r)\rangle \,
\langle t^{2l_2}(X_{(2)})\phi^{2l_2}_{l_2}(p),
\phi^{2l_2}_{l_2}(r)\rangle \\ &=
\langle t^{2l_1}\otimes t^{2l_2}(\De(X))
\phi^{2l_1,2l_2}_{l_1,l_1}(p), \phi^{2l_1,2l_2}_{l_1,l_1}(r)\rangle
\\ &=
\sum_{l=|l_1-l_2|}^{l_1+l_2} C^{2l_1,2l_2,2l}_{l_1,l_2,l}(p)
C^{2l_1,2l_2,2l}_{l_1,l_2,l}(r)
\langle t^{2l}(X) \phi^{2l}_l(p),
\phi^{2l}_l(r)\rangle,
\endaligned
\tag\eqname{\vgllinearisationformone}
$$
where $r>0$ is another parameter and
$\De(X)=\sum_{(X)}X_{(1)}\otimes X_{(2)}$.

The dual Hopf $\ast$-algebra $\Asu$ generated by the matrix elements
of the representations $t^N$, $N\in\Zp$, of $\Usu$, is known
in terms of generators and relations, cf.
\cite{\CharP}, \cite{\KoelAAM}, \cite{\KoorZSE}.
 Koornwinder \cite{\KoorZSE}
has given an explicit expression for the element in
$\Asu$ corresponding to the linear functionals considered in
\thetag{\vgllinearisationformone};
$$
\langle t^{2l}(X)\phi^{2l}_l(p),\phi^{2l}_l(r)\rangle
= {{q^{-l}}\over{(q^{2l+2};q^2)_l}}\,
\langle X, p_l(\rho; q\sqrt{{p\over r}},q\sqrt{{r\over p}},
{{-q}\over{\sqrt{pr}}}, -q\sqrt{pr}|q^2)\rangle,
$$
where $\rho\in\Asu$ is some fixed simple element, which is, up to an
affine scaling, the linear functional
$X\mapsto  \langle t^2(X)\phi^2_1(p),\phi^2_1(r)\rangle$, and
the last $\langle\cdot,\cdot\rangle$ denotes the duality
between $\Usu$ and $\Asu$. Since $\Asu$ is the dual Hopf
$\ast$-algebra, the left hand side of
\thetag{\vgllinearisationformone} corresponds
to the multiplication of the two linear functionals.
So \thetag{\vgllinearisationformone} leads to the
following identity in $\Asu$;
$$
p_{l_1}(\rho)\, p_{l_2}(\rho) =
\sum_{l=|l_1-l_2|}^{l_1+l_2} q^{l_1+l_2-l}
{{(q^{2l_1+2};q^2)_{l_1}(q^{2l_2+2};q^2)_{l_2}}\over
{(q^{2l+2};q^2)_l}} C^{2l_1,2l_2,2l}_{l_1,l_2,l}(p)
C^{2l_1,2l_2,2l}_{l_1,l_2,l}(r)\, p_l(\rho)
$$
with $p_l(\cdot)=p_l(\cdot ; q\sqrt{{p\over r}},q\sqrt{{r\over p}},
{{-q}\over{\sqrt{pr}}}, -q\sqrt{pr}|q^2)$.

The only information on $\Asu$ needed is the existence of a family
of one-dimensional representations sending $\rho$ to $\cos\th$.
Thus, applying the one-dimensional representations of
$\Asu$ and using Proposition
\thmref{\propeigvectDeXpeneNn}
proves the following linearisation coefficient formula.

\proclaim{Theorem \theoremname{\thmlinearisationcoeffAWpols}}
Let $p_l(x)= p_l(x;q\sqrt{{p\over r}},q\sqrt{{r\over p}},
{{-q}\over{\sqrt{pr}}}, -q\sqrt{pr}|q^2)$, $p,r>0$,
be defined in terms of Askey-Wilson polynomials
\thetag{\vgldefAskeyWilsonpols}. Then the coefficients
in the linearisation formula
$$
p_{l_1}(x)\, p_{l_2}(x) =
\sum_{j=0}^{2\min(l_1,l_2)} c_j \,p_{l_1+l_2-j}(x)
$$
are given by a product of two
balanced terminating ${}_4\vp_3$-series;
$$
\align
c_j = & q^{-j(j-1)}q^{j+4jl_1} (q^{2l_1+2};q^2)_{l_1-j}
(q^{2l_2+2};q^2)_{l_2} \left[ {{2l_2}\atop j}\right]_{q^2}
\\ &\times
{{ (q^{2l_1+2l_2};q^{-2})_j}\over{
(q^{2l_1+2l_2+2};q^2)_{l_1+l_2-j}}}
{{1-q^{4l_1+4l_2-4j+2}}\over{1-q^{4l_1+4l_2-2j+2}}} \\
 &\times p^{j/2}(-p^{-1}q^{-2l_1-2l_2};q^2)_j
\ {}_4\vp_3\left( {{q^{-2j},q^{-2l_2}, q^{2j-2-4l_1-4l_2},
-p^{-1}q^{-2l_2}} \atop{q^{-4l_2}, q^{-2l_1-2l_2},
-p^{-1}q^{-2l_1-2l_2}}}; q^2,q^2\right) \\
&\times
r^{j/2}(-r^{-1}q^{-2l_1-2l_2};q^2)_j
\ {}_4\vp_3\left( {{q^{-2j},q^{-2l_2}, q^{2j-2-4l_1-4l_2},
-r^{-1}q^{-2l_2}} \atop{q^{-4l_2}, q^{-2l_1-2l_2},
-r^{-1}q^{-2l_1-2l_2}}}; q^2,q^2\right).
\endalign
$$
\endproclaim

\demo{Remark \theoremname{\remthmlinearisationcoeffAWpols}}
(i) In particular, for $p=r$ the linearisation coefficients
are positive. This can already be observed without the explicit
knowledge of the linearisation coefficients, see
\cite{\KoelAAM, \S 8.3}, \cite{\KoorProb, \S 7}.

(ii) For $p=r=1$ the Askey-Wilson polynomials
$p_l(x;q,q,-q,-q|q^2)$ are the continuous $q$-Legendre polynomials
$C_l(x;q^2|q^4)$, see \cite{\AskeW, \S 4}. This is a special case
of the continuous $q$-ultraspherical polynomials introduced by
Rogers at the end of the 18th century. Rogers calculated the
linearisation coefficients for the continuous $q$-ultraspherical
polynomials, see e.g.
\cite{\AskeW, \S 4}, \cite{\GaspR, \S 8.5}, and we can go from
Theorem \thmref{\thmlinearisationcoeffAWpols} to the special case
of Rogers's result by using Andrew's summation formula, see
\cite{\GaspR, (II.17)}.
\enddemo


\enddocument